\documentclass[aap,MSNbibl,seceqn,noautosecdot,dvips]{arximspdf}

\doi{10.1214/12-AAP886} 
\volume{23}
\issue{5}
\pubyear{2013}
\firstpage{1778}
\lastpage{1816}

\makeatletter

\newcommand{\rright}{\right}
\newcommand{\lleft}{\left}
\newcommand{\rrvert}{\vert}
\newcommand{\llvert}{\vert}

\newtheorem{theorem}{Theorem}[section]
\newtheorem{corollary}[theorem]{Corollary}
\newtheorem{lemma}[theorem]{Lemma}
\newtheorem{proposition}[theorem]{Proposition}

\newproclaim{example}[theorem]{Example}
\newproclaim{definition}[theorem]{Definition}
\newproclaim{Remark}[theorem]{Remark}

\def\cal{\mathcal}

\newcommand{\EE}{{\mathbf E}}
\newcommand{\RR}{{\mathbf R}}
\newcommand{\PP}{{\mathbf P}}

\newcommand{\ww}{{\mathbf w}}
\newcommand{\pt}{{\tilde{p}}}
\newcommand{\zeroh}{{\hat{0}}}
\newcommand{\oneh}{{\hat{1}}}
\newcommand{\Kt}{{\tilde{K}}}
\newcommand{\Fc}{{\cal F}}
\newcommand{\Kc}{{\cal K}}
\newcommand{\MMc}{{\cal M}}
\newcommand{\Xc}{{\cal X}}

\newcommand{\pih}{\hat{\pi}}
\newcommand{\zh}{\hat{0}}

\newcommand{\gl}{\lambda}

\setattribute{abstract}   {width}  {287pt}

\makeatother

\begin{document}
\begin{frontmatter}

\title{Comparison inequalities and fastest-mixing Markov~chains}
\runtitle{Fastest-mixing Markov chains}

\begin{aug}
\author[A]{\fnms{James Allen} \snm{Fill}\corref{}\thanksref{t1}\ead[label=e1]{jimfill@jhu.edu}}
\and
\author[B]{\fnms{Jonas} \snm{Kahn}\ead[label=e2]{jonas.kahn@math.univ-lille1.fr}}
\runauthor{J. A. Fill and J. Kahn}
\affiliation{Johns Hopkins University and Universit\'{e} de Lille 1, CNRS}
\address[A]{Department of Applied Mathematics\\
\quad and Statistics\\
Johns Hopkins University\\
3400 N. Charles Street\\
Baltimore, Maryland 21218-2682\\
USA\\
\printead{e1}} 
\address[B]{Laboratoire Paul Painlev\'{e} (UMR 8524)\\
Universit\'{e} de Lille 1, CNRS\\
Cit\'{e} Scientifique---B\^at. M2\\
59655 Villeneuve d'Ascq Cedex\\
France\\
\printead{e2}}
\end{aug}

\thankstext{t1}{Supported by the Acheson J. Duncan Fund for the
Advancement of Research in Statistics.}

\received{\smonth{9} \syear{2011}}
\revised{\smonth{3} \syear{2012}}

%
\begin{abstract}
We introduce a new partial order on the class of stochastically
monotone Markov kernels having a given stationary distribution $\pi$ on
a given finite partially ordered state space $\Xc$. When $K \preceq L$
in this partial order we say that $K$ and $L$ satisfy a
\textit{comparison inequality}. We establish that if $K_1,\ldots, K_t$
and $L_1,\ldots, L_t$ are reversible and $K_s \preceq L_s$ for $s =
1,\ldots, t$, then $K_1 \cdots K_t \preceq L_1 \cdots L_t$. In
particular, in the time-homogeneous case we have $K^t \preceq L^t$ for
every $t$ if $K$ and $L$ are reversible and $K \preceq L$, and using
this we show that (for suitable common initial distributions) the
Markov chain $Y$ with kernel $K$ mixes faster than the chain $Z$ with
kernel $L$, in the strong sense that \textit{at every time $t$} the
discrepancy---measured by total variation distance or separation or
$L^2$-distance---between the law of $Y_t$ and $\pi$ is smaller than
that between the law of $Z_t$ and $\pi$.

Using comparison inequalities together with specialized arguments to
remove the stochastic monotonicity restriction, we answer a question of
Persi Diaconis by showing that, among all symmetric birth-and-death
kernels on the path $\Xc= \{0,\ldots, n\}$, the one (we call it the
\textit{uniform chain}) that produces fastest convergence from initial
state $0$ to the uniform distribution has transition probability $1/2$
in each direction along each edge of the path, with holding probability
$1/2$ at each endpoint.

We also use comparison inequalities:\vspace*{4pt}

\hspace*{2.5pt}(i) to identify, when $\pi$ is a given log-concave
distribution on
the path, the fastest-mixing stochastically monotone birth-and-death
chain started at $0$, and

(ii) to recover and extend a result of Peres
and Winkler that extra updates do not delay mixing for monotone spin
systems.\vspace*{4pt}

\noindent
Among the fastest-mixing chains in (i), we show that the
chain for
uniform $\pi$ is slowest in the sense of maximizing separation at every
time.
\end{abstract}

%

\end{frontmatter}

\section{\texorpdfstring{Introduction and summary.}{Introduction and summary}} \label{Sintrosummary}

A series of papers~\cite{MR2124681,MR2278445,MR2202924,MR2515797} by
Boyd, Diaconis, Xiao and coauthors
considers the following ``fastest-mixing Markov chain'' problem. A
finite graph $G = (V, E)$ is given, together with a probability
distribution $\pi$ on $V$ such that $\pi(i) > 0$ for every $i$; the
goal is to find the fastest-mixing reversible Markov chain (FMMC) with
stationary distribution $\pi$ and transitions allowed only along the
edges in $E$. This is a very important problem because of the use of
Markov chains in Markov chain Monte Carlo (MCMC), where the goal is to
sample (at least approximately) from $\pi$ and the Markov chain is
constructed only to facilitate generation of such observations as
efficiently as possible. As their criterion for FMMC, the authors
minimize SLEM (second-largest eigenvalue in modulus---sometimes also
called the absolute value of the ``largest small eigenvalue''---defined
as the absolute value of the eigenvalue of the one-step kernel with
largest absolute value strictly less than $1$), and they find the FMMC
using semidefinite programming. (More precisely, the authors of \cite
{MR2124681,MR2202924,MR2515797} do this; the author of~\cite
{MR2278445} similarly deals with continuous-time chains and minimizes
relaxation time. See these papers for further references; in
particular, related work is found in~\cite{MR2198603}.)

While most of the results in the series are numerical,
both~\cite{MR2202924} and~\cite{MR2515797} contain analytical results.
For the problem treated in~\cite{MR2202924} (which, as explained there,
has an application to load balancing for a network of processors \cite
{Diekmann97engineeringdiffusive}), the graph $G$ is a path (say, on $V
= \{0,\ldots, n\}$, with an edge joining each consecutively-numbered
pair of vertices) with a self-loop at each vertex, $\pi$ is the uniform
distribution, and it is proved that the FMMC has transition probability
$p(i, i + 1) = p(i + 1, i) = 1 / 2$ along each edge and $p(i, i) \equiv
0$ except that $p(0, 0) = 1 / 2 = p(n, n)$. [We will call this the
\textit{uniform} chain $U = (U_t)_{t = 0, 1, \ldots}$.]

The mixing time of a Markov chain can indeed be bounded using the SLEM,
which provides the asymptotic exponential rate of convergence to
stationarity. (See, e.g.,~\cite{AldousFill} for background and
standard Markov chain terminology used in this paper.) But the SLEM
provides only a surrogate for true measures of discrepancy from
stationarity, such as the standard total variation (TV) distance,
separation (sep) and $L^2$-distance. For the path problem, for example,
Diaconis (personal communication) has wondered whether the uniform
chain might in fact minimize such distances after any given number of
steps (when, for definiteness, all chains considered must start
at $0$). In this paper we show that this is indeed the case: the
uniform chain is truly fastest-mixing in a wide variety of senses.
Consider any
$t \geq0$. What we show, precisely, is that, for any birth-and-death
chain\setcounter{footnote}{1}\footnote{Arbitrary holding is allowed at each state.} $X$ having
symmetric transition kernel on the path and initial state $0$, the
probability mass function (p.m.f.) $\pi_t$ of $X_t$ majorizes the
p.m.f.
$\sigma_t$ of $U_t$.
(A definitive reference on the theory of majorization is~\cite{MR81b00002}.)
We will show using this that four examples of discrepancy from
uniformity that are larger for $X_t$ than for $U_t$ are
(i) $L^p(\pi)$-distance for any $1 \leq p \leq\infty$ (including the
standard TV and $L^2$ distances); (ii)~separation;
(iii) Hellinger distance; and (iv) Kullback--Leibler divergence.

The technique we use to prove that $\pi_t$ majorizes $\sigma_t$ is new
and remarkably simple, yet quite general. In Section~\ref{SCI} we describe
our method of comparison inequalities. We
show (Corollary~\ref{Chomog}) that if two Markov semigroups satisfy a certain
\textit{comparison inequality} at time $1$, then they satisfy the same
comparison inequality at all times $t$. We also show, in Section~\ref{Smaj}
(see especially Corollary~\ref{CTVsepL2}), how the comparison
inequality can be
used to compare mixing times---in a variety of senses---for the chains
with the given semigroups.

In Section~\ref{Sunif} we show that, in the context of the above
\textit{path-problem} (of finding the FMMC on a path), if one restricts either
(i) to monotone chains, or (ii)~to even times,
then the uniform chain satisfies a favorable comparison inequality in
comparison with any other chain in the class considered.
Somewhat delicate arguments (needed except in the case of
$L^2$-distance) specific to the path-problem allow us to remove the
parity restriction from the conclusion that the uniform chain is
fastest; see Theorem~\ref{Tpathmain}. Further, comparisons
between chains---even time-inhomogeneous ones---other than the
fastest $U$ can be carried out with our method by limiting attention
either to monotone kernels or to two-step kernels.
Indeed, our Proposition~\ref{PCI} rather generally provides a new tool
for the notoriously difficult analysis of time-inhomogeneous chains,
whose nascent quantitative theory has been advanced impressively in
recent work of Saloff-Coste and
Z\'{u}\~{n}iga~\cite{MR2340874,MR2519527,SCZmergingII,SCZwave}.

In Section~\ref{SBD} (see Theorem~\ref{TFMBD}), we generalize our
path-problem
result as follows.
Let $\pi$ be a log-concave p.m.f. on $\Xc= \{0,\ldots, n\}$. Among all
\textit{monotone} birth-and-death kernels $K$, the fastest to mix (again,
in a variety of senses)
is $K_\pi$
with (death, hold, birth) probabilities given by
\[
q_i = \frac{\pi_{i - 1}}{\pi_{i - 1} + \pi_i},\qquad r_i = \frac
{\pi^2_i - \pi_{i - 1} \pi_{i + 1}}{(\pi_{i - 1} + \pi_i) (\pi_i + \pi
_{i +
1})},\qquad
p_i = \frac{\pi_{i + 1}}{\pi_i + \pi_{i + 1}}.
\]
(This reduces to the uniform chain when $\pi$ is uniform.)

In Section~\ref{Smixing} we revisit the birth-and-death problems of
Sections~\ref{Sunif}--\ref{SBD} in terms of an alternative notion of
mixing time employed by Lov\'{a}sz and Winkler~\cite{LovaszWinkler2}.
Consider, for example, the
path-problem of Section~\ref{Sunif}. For every even value of $n$ the uniform
chain is fastest-mixing in their sense, too. But, perhaps somewhat
surprisingly, for every odd value of $n$ the uniform chain is
\textit{not} fastest-mixing in their sense; we identify the chain that is.

In Section~\ref{Sladder} we discuss a simple ``ladder'' game, where
the class
of kernels is a certain subclass of the symmetric birth-and-death
kernels considered in Section~\ref{Sunif}.

In Section~\ref{SPW} we show how comparison inequalities can recover and
extend (among other ways, to certain card-shuffling chains) a
Peres--Winkler result about slowing down mixing by skipping
(``censoring'') updates of monotone spin systems. (This is an example
of comparison inequalities applied to time-inhomogeneous chains.)

\section{\texorpdfstring{Comparison inequalities.}{Comparison inequalities}}\label{SCI}

In this section we introduce our new concept of \textit{comparison inequalities}.
Consider a p.m.f. $\pi> 0$ on a given finite partially ordered state
space $\Xc$.
We utilize the usual $L^2(\pi)$ inner product
%
%
\begin{equation}
\label{inner} \langle f, g \rangle\equiv\langle f, g \rangle_{\pi}:=
\sum_{i
\in
\Xc} \pi(i) f(i) g(i);
\end{equation}
if a matrix $K$ is regarded in the usual fashion as an operator on
$L^2(\pi)$ by regarding functions on $\Xc$ as column vectors, then the
$L^2(\pi)$-adjoint of~$K$ (also known as the time-reversal of $K$,
when $K$ is a Markov kernel) is $K^*$ with
$K^*(i, j) \equiv\pi(j) K(j, i) / \pi(i)$. Reversibility with respect
to $\pi$ for a Markov kernel $K$ is simply the condition that $K$ is
self-adjoint.

Let $\Kc$, $\MMc$ and $\Fc$ denote the respective classes of (i) Markov
kernels on $\Xc$ with stationary distribution $\pi$, (ii) nonnegative
nonincreasing functions on $\Xc$ and (iii) kernels $K$ from $\Kc$ that
are stochastically monotone (meaning that $K f \in\MMc$ for every $f
\in\MMc$).
Note for future reference that the identity kernel $I$ always belongs
to $\Fc$,
regardless of $\pi$.
Define a \textit{comparison inequality} relation $\preceq$
on $\Kc$ by declaring that $K \preceq L$ if $\langle K f, g \rangle
\leq\langle L f, g \rangle$ for every $f, g \in\MMc$, and observe
that $K \preceq L$ if and only if the time-reversals
$K^*$ and $L^*$ satisfy $K^* \preceq L^*$.
%
%
\begin{Remark}
\label{Rreduce}
(a) Clearly:
\begin{enumerate}[(iii)\hspace*{4pt}(i)]
\item[(i)] to verify a comparison inequality $K \preceq L$ by establishing
$\langle K f, g \rangle\leq\langle L f, g \rangle$, it is sufficient
to take $f$ and $g$ to be indicator functions of down-sets (i.e.,
sets $D$ such that $y \in D$ and $x \leq y$ implies $x \in D$) in the
partial order; and
\item[(ii)] if a comparison inequality holds, then the condition that $f$
and $g$ be nonnegative can be dropped, if desired.
\end{enumerate}

\begin{enumerate}[(b)]
\item[(b)] There
is an important existing notion of \textit{stochastic ordering} for
Markov kernels
on $\Xc$: we say that
$L \leq_{\mathrm{st}} K$ if $K f \leq L f$ entrywise for all
$f \in\MMc$. It is clear that
$L \leq_{\mathrm{st}} K$ implies $K \preceq L$ when $K$
and $L$ belong to $\Fc$.
But in all the examples in this paper where we prove a comparison
inequality, we \textit{do not} have stochastic ordering. This will
typically be the case for interesting examples, since the requirement
for distinct $K, L \in\Fc$ to have the same stationary distribution
makes it difficult (though not impossible) to have $L \leq_{\mathrm
{st}} K$.
\end{enumerate}
\end{Remark}
%
%
\begin{Remark}
\label{RCIpo}
The relation $\preceq$ defines a partial order on $\Kc$. Indeed,
reflexivity and transitivity are immediate, and antisymmetry follows\vadjust{\goodbreak}
because one can build a basis for functions on $\Xc$ from elements $f$
of $\MMc$, namely, the indicators of principal down-sets (i.e.,
down-sets of the form $\langle x \rangle:= \{y\dvtx  y \leq x\}$ with $x
\in\Xc$).
A proof from first principles is easy.\footnote{We need only show that
the indicator function ${\mathbf1}_{\{x\}}$
of any singleton $\{x\}$ can be written as a linear combination of
indicator functions of principal down-sets. But this can be done
recursively by starting with minimal elements $x$ and then using the identity
\[
{\mathbf1}_{\{x\}} = {\mathbf1}_{\langle x \rangle} - \sum
_{y < x} {\mathbf1}_{\{
y\}},\qquad x \in\Xc.
\]}%
\end{Remark}

We list next a few basic properties of the comparison relation $\preceq
$ on $\Kc$, showing that the relation is preserved under passages to
limits, mixtures, and direct sums. The proofs are all very easy. Note
also that the
class $\Fc$ of stochastically monotone kernels with stationary
distribution $\pi$ is closed under passages to limits and mixtures, and
also under (finite) products, but not under general direct sums as in
part~(c).
%
%
\begin{proposition}
\label{PCIbasic}
\textup{(a)} If $K_t \preceq L_t$ for every $t$ and $K_t \to K$ and $L_t
\to L$, then $K \preceq L$.

\textup{(b)} If $K_t \preceq L_t$ for $t = 0, 1$ and $0 \leq\gl\leq1$, then
\[
(1 - \gl) K_0 + \gl K_1 \preceq(1 - \gl)
L_0 + \gl L_1.
\]

\textup{(c)} Partition $\Xc$ arbitrarily into subsets $\Xc_0$ and $\Xc_1$,
and let each $\Xc_i$ inherit its partial order and stationary
distribution from $\Xc$. For $i = 0, 1$, suppose $K_i \preceq L_i$ on
$\Xc_i$. Define the kernel $K$
(resp., $L$) as the direct sum of $K_0$ and $K_1$ (resp., $L_0$
and~$L_1$). Then $K \preceq L$.
\end{proposition}

The following proposition, showing that $\preceq$ is preserved under
product for stochastically monotone reversible kernels, is the main
result of this section.

%
\begin{proposition}[(Comparison inequalities)]
\label{PCI}
Let $K_1,\ldots, K_t$ and $L_1,\ldots,\break L_t$ be reversible [i.e.,
$L^2(\pi
)$-self-adjoint] kernels all belonging to $\Fc$, and suppose that $K_s
\preceq L_s$ for $s = 1,\ldots, t$. Then the product kernels $K_1
\cdots K_t$ and $L_1 \cdots L_t$ (and their time-reversals)
belong to $\Fc$, and $K_1 \cdots K_t \preceq L_1 \cdots L_t$.
\end{proposition}

The application to time-homogeneous chains is the following immediate corollary.
%
%
\begin{corollary}
\label{Chomog}
If $K, L \in\Fc$ are
reversible and $K \preceq L$, then for every $t$ we have $K^t, L^t \in
\Fc$ and $K^t \preceq L^t$.
\end{corollary}
%

\begin{Remark}
\label{Rapplicability}
As we shall see from examples, the applicability of our new technique
of comparison inequalities is limited (i) by the monotonicity
requirement for membership in\vadjust{\goodbreak} $\Fc$ and (ii) by the extent to
which $\Fc
$ is ordered by $\preceq$. But restriction (i) in the choice of kernel
has the payoff (among others) that the perfect simulation algorithms
(see~\cite{MR1630416} for background) Coupling From The Past \cite
{MR99k60176,MR99d60081,MR99g60116,MR2001h65014} and FMMR
(Fill--Machida--Murdoch--Rosenthal)~\cite{MR99g60113,MR2001m60164}
can often be run efficiently for monotone chains.
Restriction (ii) needs to be explored thoroughly for interesting and important
examples. This paper treats a few examples, in Sections~\ref{Sunif}
(especially Section~\ref{SSmono}),~\ref{SBD}
and~\ref{SPW}.
For discussion about the relation between our comparison-inequalities
technique and existing techniques for comparing mixing times of Markov
chains, see Remark~\ref{Rcompcomp} below.
\end{Remark}

The remainder of this section is devoted to the proof of
Proposition~\ref{PCI},
which we will derive as a consequence of an extremely simple, but---as
far as we know---new, matrix-theoretic result, Proposition~\ref{Pmatrix}.

The general setting is this. We are given a positive vector $\pi\in
\RR^n$ and define the $L^2(\pi)$ inner product as at (\ref{inner}).
We are also given a set (not necessarily a subspace) $W \subseteq\RR
^n$. Let $M_n(\RR)$ denote the collection of $n$-by-$n$ real matrices. Define
\[
\Fc:= \bigl\{\mbox{matrices $A \in M_n(\RR)$ for which $W$ is
invariant}\bigr\}.
\]
(This of course means that a real matrix $A$ belongs to $\Fc$ if and
only if $A w \in W$ for every $w \in W$.)
Define a (clearly reflexive and transitive) relation $\preceq$ on
$M_n(\RR)$ by declaring that $A \preceq B$ if
\[
\langle A x, y \rangle\leq\langle B x, y \rangle\qquad\mbox{for every $x, y
\in W$}.
\]
We observe in passing (i) that $A \preceq B$ if and only if $A^*
\preceq B^*$ and (ii) that the relation $\preceq$ may fail to be
antisymmetric (but this will present no difficulty).

%
\begin{proposition}
\label{Pmatrix}
Let $A_1, A_2, B_1, B_2 \in M_n(\RR)$. Suppose that $A_2$ and $B_1^*$
both belong to $\Fc$. If $A_1 \preceq B_1$ and $A_2 \preceq B_2$, then
$A_1 A_2 \preceq B_1 B_2$.
\end{proposition}
\begin{pf}
Given $x, y \in W$, we observe
\begin{eqnarray*}
\langle A_1 A_2 x, y \rangle&\leq&\langle
B_1 A_2 x, y \rangle\qquad\mbox{because $A_2 x, y
\in W$ and $A_1 \preceq B_1$}
\\
&=& \bigl\langle A_2 x, B_1^* y \bigr\rangle
\\
&\leq&\bigl\langle B_2 x, B_1^* y \bigr\rangle
\qquad\mbox{because $x, B_1^* y \in W$ and $A_2 \preceq
B_2$}
\\
&=& \langle B_1 B_2 x, y \rangle
\end{eqnarray*}
as desired.
\end{pf}

The third (Corollary~\ref{Cc3}) of the following four easy
corollaries of
Proposition~\ref{Pmatrix} implies Proposition~\ref{PCI}
immediately, by setting $W = \MMc$ and
observing that the set of Markov kernels with stationary distribution
$\pi> 0$ is closed under both multiplication and adjoint. (Similarly,
Corollary~\ref{Chomog} is a special case of Corollary~\ref{Cc4}.)
%
%
\begin{corollary}
\label{Cc1}
Let $A_1, A_2, B_1, B_2$ be matrices all belonging to $\Fc$ with
adjoints all belonging to $\Fc$, and suppose that $A_1 \preceq B_1$ and
$A_2 \preceq B_2$. Then the matrices $A_1 A_2$ and $B_1 B_2$ and their
adjoints all belong to $\Fc$, and $A_1 A_2 \preceq B_1 B_2$.
\end{corollary}
\begin{pf}
This is immediate from the definition of $\Fc$ and
Proposition~\ref{Pmatrix}.
\end{pf}
%
%
\begin{corollary}
\label{Cc2}
Let $A_1,\ldots, A_t$ and $B_1,\ldots, B_t$ be matrices all belonging
to $\Fc$ with adjoints all belonging to $\Fc$, and suppose that $A_s
\preceq B_s$ for $s = 1,\ldots, t$. Then the matrices $A_1 \cdots A_t$
and $B_1 \cdots B_t$ and their adjoints all belong to $\Fc$, and $A_1
\cdots A_t \preceq B_1 \cdots B_t$.
\end{corollary}
\begin{pf}
This follows by induction from Corollary~\ref{Cc1}.
\end{pf}
%
%
\begin{corollary}
\label{Cc3}
Let $A_1,\ldots, A_t$ and $B_1,\ldots, B_t$ be self-adjoint matrices
all belonging to $\Fc$, and suppose that $A_s \preceq B_s$ for $s =
1,\ldots, t$. Then the matrices $A_1 \cdots A_t$ and $B_1 \cdots B_t$ (and
their adjoints) belong to $\Fc$, and $A_1 \cdots A_t \preceq B_1
\cdots B_t$.
\end{corollary}
\begin{pf}
This is immediate from Corollary~\ref{Cc2}.
\end{pf}
%
%
\begin{corollary}
\label{Cc4}
Let $A$ and $B$ be self-adjoint matrices both belonging to~$\Fc$, and
suppose that $A \preceq B$. Then, for every $t = 0, 1, 2, \ldots\,$, the
matrices $A^t$ and $B^t$ (are self-adjoint and) belong to $\Fc$ and
$A^t \preceq B^t$.
\end{corollary}
\begin{pf}
This is immediate from Corollary~\ref{Cc3} by taking $A_s \equiv A$
and $B_s
\equiv B$.
\end{pf}

\section{\texorpdfstring{Consequences of the comparison inequality, some via
majorization.}{Consequences of the comparison inequality, some via
majorization}}\label{Smaj}
In this section we focus on time-homogeneous chains and show how
comparison inequalities can be used to compare mixing times---in a
variety of senses---for chains with the given semigroups. As we shall
see in Section~\ref{SSmaj}, a useful tool in moving from a comparison
inequality to a comparison of mixing times will be the use of basic
results from the theory of majorization.

\subsection{\texorpdfstring{Comparison inequalities and domination.}{Comparison inequalities and domination}}\label{SSCIdom}
Recall from Section~\ref{SCI} that $\Fc$ denotes the class of stochastically
monotone Markov kernels on a given finite partially ordered state
space $\Xc$ that have a given $\pi$ as stationary distribution. Our
next result (Proposition~\ref{PCIdom}) gives conditions implying that
if a
comparison inequality holds between reversible kernels $K, L \in\Fc$,
then the univariate distributions of the corresponding Markov chains
satisfy corresponding stochastic inequalities. The proposition utilizes
the following definition.
%
%
\begin{definition}
\label{Ddom}
Let $(Y_t)$ and $(Z_t)$ be stochastic processes with the same finite
partially ordered state space. If for every $t$ we have $Y_t \geq Z_t$
stochastically, that is,
%
%
\begin{equation}
\label{down} \PP(Y_t \in D) \leq\PP(Z_t \in
D)\qquad \mbox{for every down-set $D$ in the partial order},\hspace*{-35pt}
\end{equation}
then we say that \textit{$Y$ dominates $Z$}.
\end{definition}
%
%
\begin{proposition}
\label{PCIdom}
Suppose that $K, L \in\Fc$ are reversible and satisfy $K \preceq L$.
If $Y$ and $Z$ are chains \textup{(i)} started in a common p.m.f. $\pih$ such that
$\pih/ \pi$ is nonincreasing and \textup{(ii)} having respective kernels $K$
and $L$, then $Y$ dominates~$Z$.
\end{proposition}
\begin{pf}
By Corollary~\ref{Chomog} for every $t$ we have $K^t, L^t \in\Fc$
and $K^t
\preceq L^t$. The desired result now follows easily.
\end{pf}

\subsection{\texorpdfstring{TV, separation and $L^2$-distance.}{TV, separation and $L^2$-distance}}\label{SSTVsepL2}
Domination (recall Definition~\ref{Ddom}) is quite useful for
comparing mixing
times in at least three standard senses.

If $d$ is some measure of discrepancy
from stationarity, then in the following theorem we write ``$Y$ mixes
faster in $d$ than does $Z$'' for the strong assertion that at every
time $t$ we have $d$ smaller for $Y$ than for $Z$.
%
%
\begin{corollary}
\label{CTVsepL2}
Consider (not necessarily reversible) Markov chains $Y$ and $Z$ with
common finite partially ordered state
space $\Xc$, common initial distribution $\pih$ and common stationary
distribution $\pi$.
Assume that $\pih/ \pi$ is nonincreasing.\vspace*{8pt}

\textup{(a) (total variation distance)}.
Suppose that $Y$ dominates $Z$ and that the time-reversal of $Y$ is
stochastically monotone.
Then $Y$ mixes faster in TV than does~$Z$.

\textup{(b) (separation)}.
Adopt the same hypotheses as in part \textup{(a)}. Then $Y$ mixes faster
in separation than does $Z$; equivalently, any fastest strong
stationary time for $Y$ is stochastically smaller (i.e., faster) than
any strong stationary time for $Z$.

\textup{(c) ($L^2$-distance)}.
Assume that $Y$ and $Z$ are reversible.
Suppose, moreover, that the two-step chain $(Y_{2 t})$ dominates $(Z_{2
t})$ and is stochastically monotone.
Then $Y$ mixes faster in $L^2$ than does $Z$.
\end{corollary}
\begin{pf}
All three results are simple applications of the domination
inequality~(\ref{down}) [which, in the case of part (c), is guaranteed
only for even values of $t$] or its immediate extension to expectations
of nonincreasing functions. We make the preliminary observation that
$\PP(Y_t = i) / \pi(i)$ is nonincreasing in $i$ for each $t$; indeed,
writing $K$ for the kernel of $Y$ we have
%
%
\begin{equation}
\label{MLR} \frac{\PP(Y_t = i)}{\pi(i)} = \sum_j
\frac{\pih(j) K^t(j, i)}{\pi
(i)} = \sum_j {K^*}^t(i,
j) \frac{\pih(j)}{\pi(j)},
\end{equation}
so the nonincreasingness claimed here follows from the monotonicity
assumptions about $\pih/ \pi$ and $K^*$.\vspace*{8pt}

(a) Choosing $D$ in (\ref{down}) to be the down-set $D = \{i\dvtx \PP(Y_t =
i) / \pi(i) > 1\}$ we find
\[
\mathrm{TV}_Y(t) = \PP(Y_t \in D) - \pi(D) \leq
\PP(Z_t \in D) - \pi(D) \leq\mathrm{TV}_Z(t).
\]

(b) We first observe
\[
\mathrm{sep}_Y(t) = \max_i \biggl[ 1 - \frac{\PP(Y_t = i)}{\pi(i)}
\biggr] = 1 - \frac{\PP(Y_t = x_1)}{\pi(x_1)}
\]
for some maximal element $x_1$ in $\Xc$. Therefore, choosing $D = \Xc
\setminus\{x_1\}$ we find
\begin{eqnarray*}
\mathrm{sep}_Y(t) &=& 1 - \frac{\PP(Y_t = x_1)}{\pi(x_1)}
\leq 1 -\frac{\PP(Z_t = x_1)}{\pi(x_1)} \\
&\leq&\max_i \biggl[ 1 - \frac
{\PP(Z_t = i)}{\pi(i)}
\biggr] = \mathrm{sep}_Z(t).
\end{eqnarray*}

(c) Using routine calculations suppressed here, one finds that the
squared $L^2(\pi)$-distance (of the density with respect to $\pi$) from
stationarity for $Y_t$ equals
\begin{eqnarray*}
\sum_i \pi(i) \biggl[ \frac{\PP(Y_t = i)}{\pi(i)} - 1
\biggr]^2 &=& \sum_{j'} \biggl[ \sum
_j \pih(j) K^{2 t}\bigl(j,
j'\bigr) \biggr] \frac
{\pih
(j')}{\pi(j')} - 1
\\
&=& \sum_{j'} \PP\bigl(Y_{2 t} =
j'\bigr) \frac{\pih(j')}{\pi(j')} - 1.
\end{eqnarray*}
But $\pih/ \pi$ is nonincreasing and $Y_{2 t} \geq Z_{2 t}$
stochastically; so this last expression does not exceed
\[
\sum_{j'} \PP\bigl(Z_{2 t} =
j'\bigr) \frac{\pih(j')}{\pi(j')} - 1 = \sum
_i \pi(i) \biggl[ \frac{\PP(Z_t = i)}{\pi(i)} - 1
\biggr]^2,
\]
which is the desired conclusion.
\end{pf}

We remark in passing that a very similar proof as for Corollary \ref
{CTVsepL2}(b) gives the analogous result for the measure of discrepancy
\[
\max_i \biggl[ \frac{\PP(Y_t = i)}{\pi(i)} - 1 \biggr],
\]
and so we also have the analogous result for the two-sided measure
%
%
\begin{equation}
\label{rpd} \max_i \biggl\llvert\frac{\PP(Y_t = i)}{\pi(i)} - 1 \biggr
\rrvert.
\end{equation}
%
%
\begin{Remark}[($L^2$-distance revisited)]
\label{RL2}
We have limited the statement of Corollary~\ref{CTVsepL2}(c) to reversible
chains for simplicity. The same proof shows, more generally, for
each $t$ that if (i) $K$ and $L$ are (not necessarily reversible)
kernels with common stationary distribution $\pi$, (ii) $\pih/ \pi$ is
nonincreasing, and (iii) $\pih K^t {K^*}^t \geq\pih L^t {L^*}^t$
stochastically, then the $L^2(\pi)$-distance from stationarity for
$Y_t$ does not exceed that for $Z_t$, where the chains $Y$ and $Z$ have
respective kernels $K$ and $L$ and common initial distribution $\pih$.
Assuming (i) and (ii), for the stochastic
inequality (iii) here it is sufficient that $K$ and $L$ and their
time-reversals $K^*$ and $L^*$ are all stochastically monotone and $K
\preceq L$.
\end{Remark}
%
%
\begin{Remark}[(Concerning eigenvalues)]
\label{Rcompcomp}
(a) if $K$ and $L$ are ergodic reversible kernels in $\Fc$ (with a
common stationary distribution $\pi$) and we have the comparison
inequality $K \preceq L$, then the SLEM for $K$ is no larger than the
SLEM for $L$. This follows rather easily from Proposition \ref
{PCIdom} and Corollary~\ref{CTVsepL2}(c) using the spectral
representations of the kernels and
the ample freedom in choice of the common initial distribution $\hat
{\pi
}$ such that $\hat{\pi} / \pi$ is nonincreasing. We omit further details.

(b) There are several existing standard techniques for comparing mixing
times of Markov chains, such as the celebrated eigenvalues-comparison
technique of Diaconis and Saloff-Coste~\cite{MR94i60074}, but none
give conclusions as strong as those available from combining
Proposition~\ref{PCIdom} and Corollary~\ref{CTVsepL2}. On the other
hand, comparison of
eigenvalues requires verifying far fewer assumptions than needed to
establish $K, L \in\Fc$ and a comparison inequality $K \preceq L$, so
our new technique is much less generally applicable.
\end{Remark}

\subsection{\texorpdfstring{Other distances via majorization.}{Other distances via majorization}}\label{SSmaj}
We now utilize ideas from majorization; see~\cite{MR81b00002} for
background on majorization and the concept of Schur-convexity used below.
For the reader's convenience we recall that, given two vectors $v$ and $w$
in ${\RR}^N$ (for some $N$), we say that \textit{$v$ majorizes $w$} if
(i) for each $k = 1,\ldots, N$ the sum of the
$k$ largest entries of $w$ is at least the corresponding sum for $v$,
and (ii) equality holds when $k = N$.
A function $\phi$ with domain $D \subseteq{\RR}^N$ is said to be
\textit{Schur-convex on $D$} if
$\phi(v) \geq\phi(w)$ whenever $v, w \in D$ and $v$ majorizes $w$.
Thus, given any two p.m.f.'s $\rho_1$ and $\rho_2$ on $\Xc$, if $\rho_1$
majorizes $\rho_2$,
then
for any Schur-convex function $\phi$ on the unit simplex (i.e., the
space of p.m.f.'s) we have $\phi(\rho_1) \geq\phi(\rho_2)$. Examples of
Schur-convex functions are given in Example~\ref{ESc} below; for each of
those examples, the inequality $\phi(\rho_1) \geq\phi(\rho_2)$ can be
interpreted as ``$\rho_2$ is closer
to $\pi$ than is~$\rho_1$.''\looseness=1

The next proposition describes one important case where we have
majorization and hence can extend the conclusions ``$Y$ mixes faster
in $d$ than does $Z$'' of Corollary~\ref{CTVsepL2} to other measures of
discrepancy $d$. Note the additional hypothesis, relative to
Corollary~\ref{CTVsepL2}, that $\pi$ is nonincreasing.
%
%
\begin{proposition}
\label{Pmaj}
Consider (not necessarily reversible) Markov chains $Y$ and $Z$ with
common finite partially ordered state
space $\Xc$, common initial distribution $\pih$, and common stationary
distribution $\pi$.
Suppose that both $\pi$ and $\pih/ \pi$ are nonincreasing.
Suppose, moreover, that $Y$ dominates $Z$ and that the time-reversal
of $Y$ is stochastically monotone.
Then, for all $t$, the p.m.f. $\pi_t$ of $Z_t$ majorizes the p.m.f. $\sigma
_t$ of $Y_t$.
\end{proposition}
\begin{pf}
As noted just above (\ref{MLR}), the ratio $\PP(Y_t = i) / \pi(i)$ is
nonincreasing in $i$; since $\pi(i)$ is also nonincreasing, so is the
product $\PP(Y_t = i)$. Hence for each $k \leq|\Xc|$ there is a
down-set $D_k$ such that
$\PP(Y_t \in D_k)$ equals the sum of the $k$ largest values of $\PP(Y_t
= i)$. Since $Y$ dominates $Z$,
inequality (\ref{down}) implies
that, for all $t$,
the p.m.f. $\pi_t$ of $Z_t$ majorizes the p.m.f. $\sigma_t$ of $Y_t$.
(This can be equivalently restated in language introduced in \cite
{coarseness}: $Z_t$ is coarser than $Y_t$, for all $t$.)
\end{pf}
%
%
\begin{corollary}
\label{CCImaj} Suppose that $K, L \in\Fc$ are reversible and satisfy $K
\preceq L$, and that their common stationary distribution $\pi$ is
nonincreasing. If $Y$ and $Z$ are chains \textup{(i)} started in a
common p.m.f. $\pih$ such that $\pih/ \pi$ is nonincreasing and
\textup{(ii)} having respective kernels $K$ and $L$, then, for all $t$,
the p.m.f. $\pi_t$ of $Z_t$ majorizes the p.m.f. $\sigma_t$ of $Y_t$.
\end{corollary}
\begin{pf}
The desired conclusion follows immediately upon combining
Propositions~\ref{PCIdom} and~\ref{Pmaj}.
\end{pf}
%
%
\begin{example}
\label{ESc}
In this example we show when $\pi$ is uniform in Proposition~\ref
{Pmaj} (or
Corollary~\ref{CCImaj}), then $Y$ mixes faster
than does $Z$ in more senses than TV, separation, and $L^2$.

Write $N$ for the size of the state space $\Xc$. Then each
of the following six functions is
Schur-convex
on the unit simplex in ${\RR}^N$:
\begin{eqnarray*}
\phi_1(v) &:=& \biggl[ N^{p - 1} \sum
_i \bigl\llvert v_i - N^{-1} \bigr
\rrvert^p \biggr]^{1 / p} \qquad\mbox{(for any $1 \leq p <
\infty$)},
\\
\phi_2(v) &:=& \max_i |N v_i - 1|,
\\
\phi_3(v) &:=& \max_i (1 - N v_i),
\\
\phi_4(v) &:=& \frac{1}{2} \sum_i
\bigl( v_i^{1 / 2} - N^{- 1 / 2} \bigr)^2,
\\
\phi_5(v) &:=& N^{-1} \sum_i
\ln\biggl( \frac{1 / N}{v_i} \biggr),
\\
\phi_6(v) &:=& \sum_i v_i
\ln(N v_i)
\end{eqnarray*}
in~\cite{MR81b00002}, Chapter 3, see
Sections I.1, I.1, A.2, I.1.b, D.5, and D.1, respectively. Therefore, if
$\rho_1$ majorizes $\rho_2$, then
$\rho_2$ is closer to $\pi$ than is $\rho_1$ in each of the following
six senses (where here $\pi$ is uniform
and we have written the discrepancy from $\pi$ for a generic
p.m.f. $\rho$):
\begin{longlist}
\item\textit{$L^p$-distance}
\[
\biggl[ \sum_i \pi(i) \biggl\llvert
\frac{\rho(i)}{\pi(i)} - 1 \biggr\rrvert^p \biggr]^{1 / p}
\]
for any $1 \leq p < \infty$;
\item\textit{$L^{\infty}$-distance}
\[
\max_i \biggl\llvert\frac{\rho(i)}{\pi(i)} - 1 \biggr\rrvert,
\]
also called \textit{relative pointwise distance};
\item\textit{separation}
\[
\max_i \biggl[ 1 - \frac{\rho(i)}{\pi(i)} \biggr];
\]
\item\textit{Hellinger distance}
\[
\frac{1}{2} \sum_i \pi(i) \biggl[ \sqrt{
\frac{\rho(i)}{\pi(i)}} - 1 \biggr]^2;
\]
\item the \textit{Kullback--Leibler divergence}
\[
D_{\mathrm{KL}}(\pi\| \rho) = - \sum_i \pi(i)
\ln\biggl[ \frac{\rho
(i)}{\pi(i)} \biggr];
\]
\item the \textit{Kullback--Leibler divergence}
\[
D_{\mathrm{KL}}(\rho\| \pi) = \sum_i \rho(i) \ln
\biggl[ \frac{\rho
(i)}{\pi(i)} \biggr].
\]
\end{longlist}
Of course, the $L^2$-distance considered in Corollary \ref
{CTVsepL2}(c) is the
special case $p = 2$ of example (i) here, and the TV distance of
Corollary~\ref{CTVsepL2}(a) amounts to the special case $p = 1$.
Relative pointwise
distance was also treated earlier without use of majorization at (\ref{rpd}).
\end{example}

\section{\texorpdfstring{Fastest mixing on a path.}{Fastest mixing on a path}}\label{Sunif}
We now specialize to the path-problem. Let $K$ be any symmetric
birth-and-death transition kernel on the path $\{0, 1,\ldots, n\}$, and
denote $K(i, i + 1) = K(i + 1, i)$ by $p_i$ [except that\vadjust{\goodbreak} $K(0, 0) = 1 -
p_0$ and $K(n, n) = 1 - p_{n - 1}$]; for example, when $n = 3$ we have
\[
K = \lleft[\matrix{ 1 - p_0 & p_0 & 0 & 0
\cr
p_0 & 1 - p_0 - p_1 & p_1 & 0
\cr
0 & p_1 & 1 - p_1 - p_2 &
p_2
\cr
0 & 0 & p_2 & 1 - p_2 } %
\rright].
\]
In this section we first show, in Sections
\ref{SSmono}--\ref{SSeven}, that if one restricts attention either:
\begin{longlist}
\item to monotone chains, or
\item to even times,
\end{longlist}
then the uniform chain $U$ with kernel $K_0$ where $p_i \equiv1 / 2$
satisfies a favorable comparison inequality in comparison with the
general $K$-chain, and we can apply all the results of Section~\ref{Smaj}.
Then, in Section~\ref{SSremoval}, we show that the parity restriction
in (ii)
can be removed to conclude that the uniform chain is, among all
symmetric birth-and-death chains, closest to uniformity (in several
senses) at all times.
In this section and the next we make use of the general observation
that a discrete-time birth-and-death chain with kernel $K$ on $\Xc= \{
0, 1,\ldots, n\}$ is monotone if and only if
%
%
\begin{equation}
\label{monocrit} K(i, i + 1) + K(i + 1, i) \leq1\qquad\mbox{for $i = 0,\ldots
, n -
1$}.
\end{equation}

Before we separate into the two cases (i) and (ii) for the
path-problem, let us note that if $f$ is the indicator of the down-set
$\{0, 1,\ldots, \ell\}$, then $K f$ satisfies
%
%
\begin{equation}
\label{Kfj} (K f)_j = %
\cases{ 1, &\quad if $0 \leq j \leq
\ell- 1$,
\cr
1 - p_{\ell}, &\quad if $j = \ell$,
\cr
p_{\ell}, &\quad if $j =
\ell+ 1$,
\cr
0, &\quad otherwise } %
\end{equation}
(with $p_n = 0$); hence if $g$ is the indicator of the down-set $\{0,
1,\ldots, m\}$, then
%
%
\begin{equation}
\label{Kfg} \langle K f, g \rangle= \frac{1}{n + 1} \times%
\cases{ m + 1, &\quad if $0 \leq m \leq\ell- 1$,
\cr
\ell+ 1 - p_\ell, &\quad if
$m = \ell$,
\cr
\ell+ 1, &\quad if $\ell+ 1 \leq m \leq n$. } %
\end{equation}

\subsection{\texorpdfstring{Restriction to monotone chains.}{Restriction to monotone chains}}\label{SSmono}
Applying (\ref{monocrit}),
our symmetric kernel $K$ is monotone if and only if $p_i \leq1 / 2$
for $i = 0,\ldots, n - 1$. Among all such choices, it is clear
that (\ref{Kfg}) is minimized when $K = K_0$. From Remark~\ref{Rreduce}(i)
it therefore follows that $K_0 \preceq K$ and hence
from Section~\ref{Smaj} (especially Corollary~\ref{CCImaj} and
Example~\ref{ESc})
that $K_0$ is fastest-mixing in several senses.\looseness=1
%
%
\begin{Remark}
\label{Rpathmono}
In fact, from (\ref{Kfg}) we see that
monotone symmetric birth-and-death kernels $K$ are monotonically
decreasing in the partial order $\preceq$ with respect to each $p_i$.
\end{Remark}

\subsection{\texorpdfstring{Restriction to even times.}{Restriction to even times}}\label{SSeven}
In the present setting of symmetric birth-and-death kernel, note that
our restriction (simply to ensure that $K$ is a kernel) on the values
$p_i > 0$ is that $p_i + p_{i + 1} \leq1$ for $i = 0,\ldots, n - 1$.
It is then routine to check that $K^2$ is (like $K$) reversible and
(perhaps unlike $K$) monotone. Indeed,
if $f$ is the indicator of the down-set $\{0, 1,\ldots, \ell\}$, then
$K^2 f$ satisfies
%
%
\begin{equation}
\label{K2fj} \bigl(K^2 f\bigr)_j = %
\cases{
1, &\quad if $0 \leq j \leq\ell- 2$,
\cr
1 - p_{\ell- 1} p_{\ell}, &\quad if $j
= \ell- 1$,
\cr
1 - 2 p_{\ell} + 2 p^2_{\ell} +
p_{\ell- 1} p_{\ell}, &\quad if $j = \ell$,
\cr
2 p_{\ell} - 2
p^2_{\ell} - p_{\ell} p_{\ell+ 1}, &\quad if $j =
\ell+ 1$,
\cr
p_{\ell} p_{\ell+ 1}, &\quad if $j = \ell+ 2$,
\cr
0, &\quad
otherwise, } %
\end{equation}
which is easily checked to be nonincreasing in $j$.

Suppose now that $g$ is the indicator of the down-set $\{0, 1,\ldots,
m\}$. Then
using (\ref{K2fj}) we can calculate, and subsequently minimize over the
allowable choices of $p_0,\ldots, p_{n - 1}$, the quantity $\langle K^2
f, g \rangle$ by considering three cases:\vspace*{8pt}

(a) Suppose $m = \ell$. Then
\[
(n + 1) \bigl\langle K^2 f, g \bigr\rangle= \ell+ (1 -
p_{\ell})^2 + p_{\ell}^2
\]
is minimized (regardless of value $\ell$) when $p_i = 1 / 2$ for $i =
0,\ldots, n - 1$.

(b) Suppose $\ell$ and $m$ differ by exactly $1$, say, $m = \ell+ 1$. Then
\[
(n + 1) \bigl\langle K^2 f, g \bigr\rangle= \ell+ (1 -
p_{\ell}) + p_{\ell} (1 - p_{\ell+ 1}) = \ell+ 1 -
p_{\ell} p_{\ell+ 1}
\]
is minimized (regardless of $\ell$) when $p_i = 1/2$ for $i = 0,\ldots,
n - 1$.

(c) Suppose $\ell$ and $m$ differ by at least $2$, say, $m \geq\ell+
2$. Then
\[
(n + 1) \bigl\langle K^2 f, g \bigr\rangle= \ell+ (1 -
p_{\ell}) + p_{\ell} + 0 = \ell+ 1
\]
does not depend on the choice of the vector ${\mathbf p}$.

From Remark~\ref{Rreduce}(i) it therefore follows that $K_0^2 \preceq K^2$
and hence (from Section~\ref{Smaj})
that $K_0^2$ is fastest-mixing in several senses. Specifically:
%
%
\begin{equation}
\label{even} \mbox{for all \textit{even} $t$, the p.m.f. $\pi_t$ of
$X_t$ majorizes the p.m.f. $\sigma_t$ of $U_t$},
\end{equation}
if $X$ and $U$ have respective kernels $K$ and $K_0$ and common
nonincreasing initial p.m.f. $\pih$.
Further, when we consider all symmetric birth-and-death chains started in
state $0$, it follows from Corollary~\ref{CTVsepL2}(c) that the chain with
kernel $K_0$ is fastest-mixing in $L^2$ (without the need to restrict
to even times, nor to monotone chains).
%
%
\begin{Remark}
\label{R2stepCI}
From the above calculations we see more generally that if $K$ and $\Kt$
are two symmetric birth-and-death kernels and for every $i$ we have
\[
\bigl\llvert p_i - \tfrac{1}{2} \bigr\rrvert\geq\bigl
\llvert\pt_i - \tfrac{1}{2} \bigr\rrvert\quad\mbox{and}\quad
p_i p_{i + 1} \leq\pt_i \pt_{i + 1},
\]
then $\Kt^2 \preceq K^2$.
\end{Remark}

\subsection{\texorpdfstring{Removal of parity restriction.}{Removal of parity restriction}}\label{SSremoval}

Throughout this subsection all chains are assumed to start at
state $0$, even when we do not explicitly declare so.
The main result of this section is the following theorem, which
extends (\ref{even}) to \textit{all} times $t = 0, 1, 2, \ldots$ and
therefore demonstrates (by Example~\ref{ESc}) that the uniform chain is
fastest to mix in a variety of senses.
%
%
\begin{theorem}
\label{Tpathmain} Let $X$ be a birth-and-death chain with state space
$\Xc= \{0, 1,\ldots, n\}$ and symmetric kernel, and let $U$ be the
uniform chain. Suppose that both chains start at $0$, and let $\pi_t$
(resp., $\sigma_t$) denote the probability mass function of $X_t$
(resp., $U_t$). Then
\[
\mbox{$\pi_t$ majorizes $\sigma_t$}\qquad \mbox{for all $t$}.
\]
\end{theorem}

Let $X$ have kernel $K$ as described at the outset of Section~\ref{Sunif}.
Let $\Pi_t$ and $\Sigma_t$ denote the cumulative distribution functions
(c.d.f.'s) corresponding to $\pi_t$ and $\sigma_t$, respectively: for example,
\[
\Sigma_t(j):= \sum_{i = 0}^j
\sigma_t(i) = \PP(U_t \leq j).
\]
From Section~\ref{SSeven} we already know that if $t$ is even, then
%
%
\begin{equation}
\label{cdfs} \Pi_t(i) \geq\Sigma_t(i)
\qquad\mbox{for all $i$},
\end{equation}
because then $\pi_t$ majorizes $\sigma_t$ and both p.m.f.'s are nonincreasing.

We build to the proof of Theorem~\ref{Tpathmain} by means of a
sequence of lemmas.
We start with a few results about the uniform chain.
%
%
\begin{lemma}
\label{L12}
\textup{(a)} For every time $t$, the p.m.f. $\sigma_t$ is nonincreasing on
its domain $\{0,\ldots, n\}$.

\textup{(b)} The distribution ``evolves by steps of two,'' depending on
parity: for $i = 0,\ldots, n - 1$ we have
\[
\sigma_t(i) = \sigma_t(i+1) \qquad\mbox{if $t + i$ is
odd}.
\]

\textup{(c)} For
every time $t$, the c.d.f. $\Sigma_t$ is concave (at integer arguments):
%
%
\begin{equation}
\label{concave} 2 \Sigma_t(i) \geq\Sigma_t(i + 1) +
\Sigma_t(i - 1),\qquad i \geq0.
\end{equation}

\textup{(d)} Inequality (\ref{concave}) is an equality if $i \geq0$
and $t$ and $i$ have opposite parity:
\[
2 \Sigma_t(i) = \Sigma_t(i + 1) + \Sigma_t(i
- 1) \qquad\mbox{if $t + i$ is odd}.
\]
\end{lemma}
\begin{pf}
(a) This was proved in a more general setting just above (\ref{MLR}).

(b) We use induction on $t$. The base case $t = 0$ is obvious ($0 = 0$).

Using the induction hypothesis at the second equality, we conclude,
when $t$ and $i \in\{1,\ldots, n - 1\}$ have opposite parity, that
\[
\sigma_t(i) = \tfrac{1}{2} \bigl[ \sigma_{t - 1}(i - 1)
+ \sigma_t(i + 1) \bigr] = \tfrac{1}{2} \bigl[
\sigma_{t - 1}(i) +\sigma_{t - 1}(i+2) \bigr] =
\sigma_t(i + 1).\vadjust{\goodbreak}
\]
Similarly, when $t$ is odd we have
\[
\sigma_t(0) = \tfrac{1}{2} \bigl[ \sigma_{t - 1}(0) +
\sigma_t(1) \bigr] = \tfrac{1}{2} \bigl[\sigma_{t - 1}(0)
+\sigma_{t - 1}(2) \bigr] = \sigma_t(1).
\]

(c) We first remark that it is well known that (\ref{concave}) is
indeed equivalent to concavity of $\Sigma_t$ at integer arguments. We
then need only note that (\ref{concave}) is merely a rewriting of the
monotonicity in part (a).
Indeed,
%
%
\begin{eqnarray}
\label{sigSig} 2 \Sigma_t(i) & = & \Sigma_t(i + 1) +
\Sigma_t(i - 1) + \sigma_t(i) - \sigma_t(i +
1)
\nonumber\\[-8pt]\\[-8pt]
& \geq &\Sigma_t(i + 1) + \Sigma_t(i - 1).
\nonumber
\end{eqnarray}

(d) Again using the equality at (\ref{sigSig}), this is merely a
rewriting of the ``steps of two'' evolution in part (b).
\end{pf}
%
%
\begin{lemma}
\label{L15}
For any time $t$ and any state $i$, if $\Pi_t(j) \geq\Sigma_t(j)$ for
all states $j$ in $[i - 2, i + 2]$, then $\Pi_{t + 2}(i) \geq\Sigma_{t
+ 2}(i)$.
\end{lemma}
\begin{pf}
In the following calculations, we lean heavily on the fact that we are
dealing with birth-and-death chains.
Utilizing natural notation such as $K^2(h, \leq i)$ for $\sum_{j \leq
i} K^2(h, j)$, we find using summation by parts that
\begin{eqnarray*}
\Pi_{t + 2}(i) &=& \sum_{h = 0}^{i + 2}
\pi_t(h) K^2(h, \leq i)
\\
&=& \sum_{j = 0}^{i + 2} \Pi_t(j)
\bigl[ K^2(j, \leq i) - K^2(j + 1, \leq i) \bigr]
\\
&=& \sum_{j = i - 2}^{i + 2} \Pi_t(j)
\bigl[ K^2(j, \leq i) - K^2(j + 1, \leq i) \bigr].
\end{eqnarray*}
Recalling that $K^2$ is monotone, the expression in square brackets
here is nonnegative, so first by hypothesis and then by reversing the
above steps (now with $\Sigma$ in place of $\Pi$) we have
\[
\Pi_{t + 2}(i) \geq\sum_{j = i - 2}^{i + 2}
\Sigma_t(j) \bigl[ K^2(j, \leq i) - K^2(j + 1,
\leq i) \bigr] = \sum_{h = 0}^{i + 2}
\sigma_t(h) K^2(h, \leq i).
\]
But $K^2_0 \preceq K^2$ (as noted in Section~\ref{SSeven}) and
$\sigma_t$ is
nonincreasing [Lemma~\ref{L12}(a)], so we finally conclude
\[
\Pi_{t + 2}(i) \geq\sum_{h = 0}^{i + 2}
\sigma_t(h) K_0^2(h, \leq i) =
\Sigma_{t + 2}(i)
\]
as desired.\vadjust{\goodbreak}
\end{pf}

An immediate consequence is the following:
%
%
\begin{lemma}
\label{L16}
If $p_0 \leq1 / 2$, then $\Pi_t(i) \geq\Sigma_t(i)$ for all times $t$
and all states~$i$.
\end{lemma}
\begin{pf}
As previously discussed, we need only consider odd times, for which the
proof is immediate by induction using Lemma~\ref{L15} once the basis
$t = 1$ is handled. But indeed
\[
\Pi_1(0) = 1 - p_0 \geq\tfrac{1}{2} =
\Sigma_1(0)
\]
and $\Pi_1(i) = 1 = \Sigma_1(i)$ for $i \geq1$.
\end{pf}

We can also prove that $\Pi_t(i) \geq\Sigma_t(i)$ for all $t$ if the
transition probability from $i$ to $i+1$ is sufficiently low:
%
%
\begin{lemma}
\label{L17}
For any state $i$ such that $p_i \leq1 / 2$, we have $\Pi_t(i) \geq
\Sigma_t(i)$ for all times $t$.
\end{lemma}
\begin{pf}
We begin with the observation that, by last-step analysis,
\[
\Pi_t(i) = \Pi_{t - 1}(i - 1) + \pi_{t - 1}(i) (1 -
p_i) + \pi_{t - 1}(i + 1) p_i,
\]
which can be rewritten in terms of c.d.f.'s as
\[
\Pi_t(i) = p_i \Pi_{t - 1}(i + 1) + (1 - 2
p_i) \Pi_{t - 1}(i) + p_i \Pi_{t - 1}(i -
1)
\]
in general and as
\[
\Sigma_t(i) = \tfrac{1}{2} \Sigma_{t - 1}(i + 1) +
\tfrac{1}{2} \Sigma_{t -
1}(i - 1)
\]
for the uniform chain.

Again we need only prove the lemma for odd times $t$, and then we find
\begin{eqnarray*}
\Pi_t(i) & = & p_i \Pi_{t - 1}(i + 1) + (1 - 2
p_i) \Pi_{t - 1}(i) + p_i \Pi_{t -
1}(i -
1)
\\
& \geq & p_i \Sigma_{t - 1}(i + 1) + (1 - 2 p_i)
\Sigma_{t - 1}(i) + p_i \Sigma_{t - 1}(i - 1)
\\
& \geq &\tfrac{1}{2} \Sigma_{t - 1}(i + 1) + \tfrac{1}{2}
\Sigma_{t-1}(i - 1)
\\
&= &\Sigma_t(i),
\end{eqnarray*}
where we know the first inequality holds because $t - 1$ is even
(whence $\Pi_{t - 1}$ dominates $\Sigma_{t - 1}$) and $p_i \leq1 / 2$,
and the second inequality follows from
concavity of $\Sigma_{t - 1}$ [Lemma~\ref{L12}(c)] again using $p_i
\leq1
/ 2$.
\end{pf}

We can now combine Lemmas~\ref{L15} and~\ref{L17} to prove:
%
%
\begin{lemma}
\label{L18}
If $p_i \leq1 / 2$ and $p_{i + 1} \leq1 / 2$, then for all times $t$
we have
%
%
\begin{equation}
\label{i+2} \Pi_t(j) \geq\Sigma_t(j)\qquad\mbox{for all $j
\geq i + 2$}.
\end{equation}
\end{lemma}
\begin{pf}
We need only consider odd times, and we proceed by induction on~$t$.
For $t = 1$ we have $\Pi_1(j) = 1 = \Sigma_1(j)$ for all $j \geq2$; so
we move on to the induction step.\vadjust{\goodbreak}

Suppose that (\ref{i+2}) holds with $t$ replaced by $t - 2$. Use of
Lemma~\ref{L17} then ensures that we in fact have $\Pi_{t - 2}(j)
\geq
\Sigma_{t - 2}(j)$ for all $j \geq i$. Hence for any $j \geq i + 2$ we
have $\Pi_{t - 2}(\ell) \geq\Sigma_{t - 2}(\ell)$ for all $\ell
\in[j
- 2, j + 2]$ and therefore, by Lemma~\ref{L15}, $\Pi_t(j) \geq
\Sigma_t(j)$.
\end{pf}
%
%
\begin{lemma}
\label{L13}
If $t + i$ is even, then
\[
\Pi_t(i) \geq\Sigma_t(i).
\]
\end{lemma}
\begin{pf}
We may assume that $t$ and $i$ are odd. In light of Lemma~\ref{L16}, we
may also assume $p_0 > 1 / 2$.
Let $2 \ell$ be the first state where the alternation of $p_i$'s
greater than and no greater than $1 / 2$ is broken:
%
%
\begin{eqnarray}
\label{break} &\displaystyle p_{2 \ell} \leq\tfrac{1}{2},&
\nonumber\\[-8pt]\\[-8pt]
&\displaystyle \forall0 \leq m < \ell\qquad  p_{2 m} > \tfrac{1}{2} \quad\mbox{and}\quad
p_{2 m + 1} \leq\tfrac{1}{2}.&
\nonumber
\end{eqnarray}
(If there is no such break, we define $2 \ell$ to be $n + 1$ or $n + 2$
according as $n$ is odd or even.) Notice that the break can happen only at
an even state, since two consecutive $p_i$'s cannot both exceed $1 / 2$.

Since $i$ is odd, we have either $i < 2 \ell$ or $i > 2 \ell$. In the
former case, condition (\ref{break}) implies $p_i \leq1 / 2$, and
Lemma~\ref{L17} proves that $\Pi_t(i) \geq\Sigma_t(i)$. In the
latter case, we
must have $2 \ell\leq n - 1$ in order for $i$ to be a state; we then
observe that $p_{2 \ell- 1} \leq1 / 2$ and $p_{2 \ell} \leq1 / 2$,
and then $\Pi_t(i) \geq\Sigma_t(i)$ by Lemma~\ref{L18}.
\end{pf}

We are now prepared to complete the proof of Theorem~\ref{Tpathmain}.
\begin{pf*}{Proof of Theorem~\ref{Tpathmain}}
Because the c.d.f. inequality (\ref{cdfs}) holds when either $t$ is even
or (by Lemma~\ref{L13}) when $t + i$ is even, we need only establish the
asserted majorization when $t$ is odd and $i$ is even. Indeed, in that
case using Lemma~\ref{L12}(d) we have
\begin{eqnarray*}
\Sigma_t(i) &=& \tfrac{1}{2} \bigl[\Sigma_t(i - 1) +
\Sigma_t(i + 1)\bigr] \leq \tfrac{1}{2} \bigl[
\Pi_t(i - 1) + \Pi_t(i + 1)\bigr]
\\
&\leq& \Pi_t(i - 1) + \max\bigl\{ \pi_t(i),
\pi_t(i + 1) \bigr\},
\end{eqnarray*}
and so there exist $i + 1$ entries of the vector $\pi_t$ whose sum is
at least $\Sigma_t(i)$. We conclude that $\pi_t$ majorizes $\sigma_t$,
as asserted.
\end{pf*}
%
%
\begin{Remark}
(a) The multiset of values $\{ \PP_i(U_t = j)\dvtx j \in\{0,\ldots, n \}
\}
$ for the uniform chain $U$ started in state $i$ does not depend on
$i \in\{0,\ldots, n\}$; therefore, the uniform chain minimizes various
distances from stationarity (including all those listed in Example~\ref{ESc})
not only when the starting state is $0$ but in the worst case over all
starting states (and indeed over all starting distributions).\vadjust{\goodbreak}

To see the asserted invariance in starting state, consider simple
symmetric random walk $V$ on the cycle $\{0,\ldots, 2n + 1\}$, with
transition probability $1 / 2$ in each direction between adjacent
states (modulo $2 n + 2$). Then for every $i, j \in\{0,\ldots, n\}$ we
have (by regarding states $n + 1,\ldots, 2n + 1$ as ``mirror
reflections'' of the states $n,\ldots, 0$, resp.)
\[
\PP_i(U_t = j) = \PP_i(V_t = j)
+ \PP_i(V_t = 2 n + 1 - j),
\]
where at most one of the two terms on the right---namely, the one with
$j - i \equiv t$ (modulo $2$)---is positive. Thus, as multisets of $2 n
+ 2$ elements each, we have the equality
\[
\bigl\{\PP_i(U_t = j)\dvtx j \in\{0,\ldots, n\}\bigr\} \cup
\{0,\ldots, 0\} = \bigl\{ \PP_i(V_t = j)\dvtx j \in\{0,\ldots,
2 n + 1\}\bigr\},
\]
where the multiset $\{0,\ldots, 0\}$ on the left here has (of course)
$n + 1$ elements. Since the multiset on the right clearly does not
depend on $i$, neither does $\{\PP_i(U_t = j)\dvtx\allowbreak j \in\{0,\ldots, n\}\}$.

(b) The SLEM (second-largest eigenvalue in modulus) is an asymptotic
measure (in the worst case over starting states) of distance from
stationarity. Accordingly, by remark (a), the uniform chain minimizes
SLEM among all symmetric birth-and-death chains. Thus we recover the
main result of~\cite{MR2202924}.
\end{Remark}

\section{\texorpdfstring{Fastest-mixing monotone birth-and-death chains.}{Fastest-mixing monotone birth-and-death chains}}\label{SBD}
Let $n$ be a positive integer and consider the state space $\Xc= \{
0,\ldots, n\}$.
Let $\pi$ be a log-concave distribution on $\Xc$, and consider the
class of discrete-time monotone birth-and-death
chains with state space $\Xc$ and stationary distribution $\pi$,
started in state $0$.
In this section we identify the fastest-mixing stochastically monotone
chain in this class as having kernel (call it
$K_{\pi}$) with (death, hold, birth) probabilities $(q_i, r_i, p_i)$
given for $i \in\Xc$ by
%
%
\begin{equation}
\label{fastpi}\qquad q_i = \frac{\pi_{i - 1}}{\pi_{i - 1} + \pi_i},\qquad r_i =
\frac
{\pi^2_i - \pi_{i - 1} \pi_{i + 1}}{(\pi_{i - 1} + \pi_i) (\pi_i + \pi
_{i +
1})},\qquad p_i = \frac{\pi_{i + 1}}{\pi_i + \pi_{i + 1}}
\end{equation}
with $\pi_{-1}:= 0$ and $\pi_{n + 1}:= 0$.
In Section~\ref{SSpifixed} we first find the FMMC when $\pi$ is held fixed;
then in Section~\ref{SSpivarying} we show that,
when $\pi$ is allowed to vary, taking it to be uniform gives the
slowest mixing in separation.

Throughout, we make heavy use of reversibility. Recall that any
irreducible birth-and-death chain on $\Xc$ is reversible with respect
to its unique stationary distribution $\pi$.

\subsection{\texorpdfstring{The FMMC when $\pi$ is fixed.}{The FMMC when $\pi$ is fixed}}\label{SSpifixed}
The main result of this subsection is the following comparison
inequality; and then Proposition~\ref{PCIdom} and Corollary~\ref
{CTVsepL2} establish
three senses (TV, separation, and $L^2$) in which the chain with kernel
$K_{\pi}$ is fastest-mixing.
%
%
\begin{theorem}
\label{TFMBD}
Let $\pi$ be log-concave on $\Xc= \{0,\ldots, n\}$.
Let $K_{\pi}$ have (death, hold, birth) probabilities\vadjust{\goodbreak} $(q_i, r_i, p_i)$
given by (\ref{fastpi}).
Then $K_{\pi}$ is a monotone birth-and-death kernel with stationary
distribution $\pi$, and
$K_{\pi} \preceq K$ for any such kernel~$K$.
\end{theorem}
\begin{pf}
Since for each $i$ the numbers $q_i, r_i, p_i$ are nonnegative ($r_i$
because of the log-concavity of $\pi$) and sum to unity, $K_{\pi}$ is
indeed a birth-and-death kernel. Since $\pi_i p_i \equiv\pi_{i + 1}
q_{i + 1}$, it is reversible with stationary distribution $\pi$. Since
$p_i + q_{i + 1} \equiv1$, it satisfies the inequality (\ref
{monocrit}) and so is monotone.

We now consider monotone birth-and-death kernels $K$ with stationary
distribution $\pi$ and \textit{general} $(q_i, r_i, p_i)$.
We prove $K_{\pi} \preceq K$ by extending the calculations in
Section~\ref{Sunif} and in particular in Section~\ref{SSmono}.
Note that if $f$ is the indicator of the down-set $\{0, 1,\ldots, \ell
\}
$, then $K f$ satisfies
%
%
\begin{equation}
\label{Kfjagain} (K f)_j = %
\cases{ 1, &\quad if $0 \leq j
\leq\ell- 1$,
\cr
1 - p_{\ell}, &\quad if $j = \ell$,
\cr
q_{\ell+ 1}, &\quad if
$j = \ell+ 1$,
\cr
0, &\quad otherwise; } %
\end{equation}
hence if $g$ is the indicator of the down-set $\{0, 1,\ldots, m\}$, then
%
%
\begin{equation}
\label{Kfgagain} \langle K f, g \rangle= %
\cases{\displaystyle  \sum
_{j = 0}^m \pi_j, &\quad if $0 \leq m \leq\ell-
1$,
\vspace*{2pt}\cr
\displaystyle \sum_{j = 0}^{\ell} \pi_j -
\pi_{\ell} p_{\ell}, &\quad if $m = \ell$,
\vspace*{2pt}\cr
\displaystyle \sum
_{j = 0}^{\ell} \pi_j, &\quad if $\ell+ 1 \leq m
\leq n$. } %
\end{equation}
Monotonicity (\ref{monocrit}) requires precisely that for each $\ell=
0,\ldots, n - 1$ we have
\[
p_{\ell} \biggl( 1 + \frac{\pi_{\ell}}{\pi_{\ell+ 1}} \biggr) = p_{\ell
} +
q_{\ell+ 1} \leq1,
\]
so clearly $K_{\pi} \preceq K$.
\end{pf}
%
%
\begin{Remark}
\label{RBDmono}
We see more generally that the kernels $K \in\Fc$ are nonincreasing
(in $\preceq$) in each $p_i$ and that
$p_i = \pi_{i + 1} / (\pi_i + \pi_{i + 1})$ maximizes $p_i$ subject to
the monotonicity constraint. (This remark generalizes Remark \ref
{Rpathmono}.) We observe in passing that the identity kernel $I$ is
the top element (i.e., unique
maximal element) in the restriction of the comparison-inequality
partial order $\preceq$ to monotone birth-and-death chains.
\end{Remark}
%
%
\begin{example}
\label{Erw}
Suppose
that the stationary p.m.f. is proportional to $\pi_i \equiv\rho^i$, that
is, is either truncated geometric (if $\rho< 1$) or its reverse (if
$\rho> 1$) or uniform\vadjust{\goodbreak} (if $\rho= 1$). Then the kernel $K_{\pi}$
corresponds to biased random walk,
%
%
\begin{equation}
\label{rw} q_i \equiv q:= \frac{1}{1 + \rho},\qquad r_i
\equiv0,\qquad p_i \equiv p:= \frac{\rho}{1 + \rho}
\end{equation}
with the endpoint exceptions, of course, that $q_0 = 0$, $r_0 = q$,
$r_n = p$, $p_n = 0$.
\end{example}

\subsection{\texorpdfstring{Slowest FMMC: The uniform chain.}{Slowest FMMC: The uniform chain}}\label{SSpivarying}
In this
subsection we consider the monotone FMMCs given by (\ref{fastpi}) for
log-concave p.m.f.'s $\pi$ and show
(Theorem~\ref{Tpivarying}) that
the uniquely slowest to mix in separation
(at every time $t$) is obtained by setting $\pi= \mbox{uniform}$.
Our first two results of this subsection consider ergodic
birth-and-death chains and their so-called strong stationary duals and
do not need any assumption about log-concavity of $\pi$. By ``ergodic''
we mean that the chain is assumed to be aperiodic, irreducible, and
positive recurrent (the third of which follows automatically from the
first two since our state space is finite) and so settles down to its
unique stationary distribution.
%
%
\begin{proposition}
\label{Pdual}
Let $X$ be an ergodic monotone birth-and-death chain on $\Xc= \{0,\ldots
, n\}$ with stationary p.m.f. $\pi$,
(death, hold, birth) transition probabilities $(q_i, r_i, p_i)$ satisfying
%
%
\begin{equation}
\label{equals1} q_{i + 1} + p_i = 1\qquad(i = 0,\ldots, n -
1)
\end{equation}
and initial state $0$.
Let $H$ denote the c.d.f. corresponding to $\pi$, with $H_{-1}:= 0$, and set
%
%
\begin{equation}
\label{dual}\qquad q^*_i= \frac{H_{i - 1}}{H_i} p_i,\qquad
r^*_i = 0,\qquad p^*_i = \frac
{H_{i + 1}}{H_i} q_{i + 1}
\qquad\mbox{($i = 0,\ldots, n - 1$)}.
\end{equation}
Then
\[
\operatorname{sep}(t) = \PP(T > t) \qquad(t = 0, 1, \ldots),
\]
where the random variable $T$ is the hitting time of state $n$ for the
birth-and-death chain $X^*$ with initial
state $0$ and transition probabilities (\ref{dual}).
\end{proposition}
\begin{pf}
The chain $X^*$ is called the strong stationary dual (SSD) of~$X$, and
the proposition is an immediate consequence of SSD theory
\cite{DF}, Section~4.3.
\end{pf}
%
%
\begin{example}
\label{Epvsq}
For a biased random walk as discussed in Example~\ref{Erw}, the dual
kernel is
\begin{eqnarray*}
q^*_i&=& \frac{1 - \rho^i}{1 - \rho^{i + 1}} \times\frac{\rho}{1 +
\rho
},\qquad
r^*_i = 0,\\
p^*_i &=& \frac{1 - \rho^{i + 2}}{1 - \rho^{i + 1}} \frac{1}{1
+ \rho}
\qquad\mbox{($i = 0,\ldots, n - 1$)}.
\end{eqnarray*}
It is easy to check that we obtain the \textit{same} dual kernel for
ratio $\rho^{-1}$ as for $\rho$. Thus if $q$ and $p$ are interchanged\vadjust{\goodbreak}
in a biased random walk with no holding except at the endpoints, then
the two chains mix equally quickly in separation.

This can be seen another way: more generally,
if the state space is a partially ordered set possessing both bottom
($\zh$) and top ($\oneh$) elements, then for any
ergodic kernel $K$ such that both $K$ and the time-reversal $\Kt$ are
stochastically monotone,
\textit{the chain $K$ from $\zh$ and the chain $\Kt$ from $\oneh$ mix
equally quickly in separation.}
Indeed, it is easy to see that for every $t$ we have, in obvious notation,
\[
\operatorname{sep}_{\zh}(t) = 1 - \frac{K^t(\zh, \oneh)}{\pi_{\oneh}} =
1 -
\frac{\Kt^t(\oneh, \zh)}{\pi_{\zh}} = \widetilde{\operatorname
{sep}}_{\oneh}(t).
\]
\end{example}
%
%
\begin{lemma}
\label{LBDsepcomparison}
Let $K$ and $L$ be two ergodic monotone birth-and-death chains on $\Xc
= \{0,\ldots, n\}$, both started at $0$, with possibly different
stationary distributions. Suppose that $K(i + 1, i) + K(i, i + 1) = 1 =
L(i + 1, i) + L(i,\allowbreak i + 1)$. Consider the notation of (\ref{dual}) and
suppose also that $p^*_i$ arising from $Y$ is at least $p^*_i$ arising
from $Z$ for all
$i = 0,\ldots, n$. Then $Y$ mixes faster in separation\footnote{Recall
our terminological convention stated in the paragraph preceding
Corollary~\ref{CTVsepL2}.} than does~$Z$.
\end{lemma}
\begin{pf}
Let $Y^*$ and $Z^*$ be the corresponding SSDs, as in Proposition~\ref{Pdual}.
An obvious coupling gives $Y^*_t \geq Z^*_t$ for every $t$, and the
lemma follows. It is worth pointing out that while the dual chains may
not be monotone, this causes no problem with the coupling because
$Y^*_t$ and $Z^*_t$ must have the same parity for every $t$; that's
because the holding probabilities for both dual chains all vanish.
\end{pf}

Next, given a FMMC for log-concave $\pi$, we show that it mixes faster
in separation than does a certain biased random walk.
%
%
\begin{theorem}
\label{TBDvsRW}
Consider the fastest-mixing monotone birth-and-death chain $X$ with
log-concave stationary p.m.f. $\pi$,
kernel (\ref{fastpi}), and initial state~$0$. Define
\[
\rho_i:= \pi_{i + 1} / \pi_i \qquad(i = 0,\ldots,
n - 1),
\]
and suppose that $i = i_0$ minimizes $| {\ln\rho_i} |$. Then $X$ mixes
faster in separation than does the biased random walk (\ref{rw}) with
$\rho$ set to
$\rho_{i_0}$.
\end{theorem}
\begin{pf}
Log-concavity is precisely the condition that $\rho_k$ is
nonincreasing in~$k$. Hence $p^*_i$ satisfies
%
%
\begin{eqnarray}\label{pstarineq}
p^*_i &=& \frac{H_{i + 1}}{H_i} \frac{\pi_i}{\pi_i + \pi_{i + 1}}
\nonumber
\\
&=& \biggl( 1 + \rho_i \frac{\pi_i}{H_i} \biggr) \times
\frac{1}{1
+ \rho_i} = \biggl( 1 + \frac{\rho_i}{\sum_{j = 0}^i \prod_{k = j}^{i - 1}
\rho_k^{-1}} \biggr) \times
\frac{1}{1 + \rho_i}
\\
&\geq& \biggl( 1 + \frac{\rho_i}{\sum_{j = 0}^i \rho
_i^{- (i - j)}} \biggr) \times
\frac{1}{1 + \rho_i} = f_i(\rho_i),\nonumber
\end{eqnarray}
where the function
%
%
\begin{equation}
\label{fidef} f_i(\rho):= \frac{1 - \rho^{i + 2}}{1 - \rho^{i + 1}}
\frac{1}{1
+ \rho}
\qquad\biggl[ \mbox{with }f_i(1):= \frac{i + 2}{2 (i + 1)} \biggr]
\end{equation}
satisfies $f_i(\rho^{-1}) \equiv f_i(\rho)$ and can be shown by
induction on $i$
to be nonincreasing in $\rho\leq1$ (and strictly so for $i \geq1$).
The induction step uses the fact that
\[
f_i(\rho) = 1 - \frac{\rho}{(1 + \rho)^2 f_{i - 1}(\rho)}
\]
together with the induction hypothesis and the (strict) increasingness
of the function $\rho\mapsto\rho/ (1 + \rho)^2$ for $\rho\leq1$. Therefore
\[
p^*_i \geq f_i(\rho_{i_0}),
\]
and this last expression is the dual birth probability from state $i$
for the biased random walk with ratio $\rho_{i_0}$. The conclusion of
the theorem now follows from Lem\-ma~\ref{LBDsepcomparison}.
\end{pf}

So the question as to which of the FMMCs (\ref{fastpi}) is slowest to
mix is reduced to finding the slowest biased random walk. But we have
already done the calculations needed to prove the following result:
%
%
\begin{theorem}
\label{TRWvsunif}
Consider biased random walks as in Example~\ref{Erw}, each with initial
state $0$. The walks are monotonically slower to mix in separation as
$\min\{p / q, q / p\}$ increases.
\end{theorem}
\begin{pf}
We have already noted at Example~\ref{Epvsq} that the speed of mixing is
invariant under interchange of $p$ and $q$. Moreover, as $\rho= p / q$
increases over $(0, 1]$, the chains are monotonically slower to mix in
separation because we have equality in (\ref{pstarineq}) and hence
\[
p^*_i = f_i(\rho),
\]
which (as shown in the proof of Theorem~\ref{TBDvsRW}) is
nonincreasing in
$\rho\leq1$.
\end{pf}

The next theorem is the main result of the subsection and is an
immediate corollary of Theorems~\ref{TBDvsRW} and~\ref{TRWvsunif}.
%
%
\begin{theorem}
\label{Tpivarying}
Among the fastest-mixing monotone birth-and-death\break chains (\ref{fastpi})
with initial state $0$ and log-concave stationary p.m.f. $\pi$, the
uniform chain is slowest to mix in separation.\vadjust{\goodbreak}
\end{theorem}
%
%
\begin{Remark}
\label{howfast}
How fast \textit{does} an ergodic monotone birth-and-death chain mix in
separation? We have addressed this question in general in
Proposition~\ref{Pdual}
and in the last sentence of Example~\ref{Epvsq}. The biased random
walk (\ref{rw})
is treated in some detail in~\cite{MR373604}, Section XVI.3. We
note:\vspace*{8pt}

(a) The eigenvalues, listed in decreasing order, are $1$ and
\[
2 \sqrt{p q} \cos\frac{\pi j}{n + 1} \qquad(j = 1,\ldots, n).
\]

(b) Fix $\rho$ and consider $n \to\infty$. Let $\mu= |p - q|$ denote
the size of the drift of the walk. If $\mu\neq0$ (i.e., $\rho\neq
1$), there is a ``cutoff phenomenon'' for separation at time $t = \mu n
+ c_{\rho} n^{1 / 2}$. This means (roughly put) that separation is
small at that time $t$ when $c_{\rho}$ is near $- \infty$ and large
when it is near $+ \infty$, with the subscript in $c_{\rho}$ indicating
that the definition of ``near'' depends on $\rho$.

(c) If $\rho= 1$ (the uniform chain), it takes time of the larger
order $n^2$ for separation to drop from near $1$ to near $0$, and in
this case there is no cutoff phenomenon.
\end{Remark}

\section{\texorpdfstring{Lov\'{a}sz--Winkler mixing times.}{Lov\'{a}sz--Winkler mixing times}}\label{Smixing}

In previous sections we have discussed mixing in terms of TV,
separation, $L^2$ and other functions measuring discrepancy. An
alternative description of speed of convergence is provided by mixing
times as defined by
Lov\'{a}sz and Winkler~\cite{LovaszWinkler2}; according to their
definition (reviewed below), and unlike for our previous notions of
mixing, one number [``the mixing time,'' $T_{\mathrm{mix}}(X)$] is assigned
to each chain $X$.

In this section we compute $T_{\mathrm{mix}}(X)$ for any irreducible
birth-and-death chain $X$ started at $0$ and then revisit the FMMC
problems of the preceding two sections using $T_{\mathrm{mix}}$ as our
criterion. One highlight is this: for the path-problem on $\Xc=
\{0,\ldots, n\}$, we show that the uniform chain is the fastest-mixing
symmetric birth-and-death chain in the sense of Lov\'{a}sz and
Winkler~\cite{LovaszWinkler2} if and only if $n$ is even, and we identify the
fastest chain when $n$ is odd.

According to the definition in~\cite{LovaszWinkler2},
the \textit{mixing time} for any irreducible (discrete-time) finite-state
Markov chain $X$ having stationary
distribution $\pi$ is the (attained) infimum of expectations of
randomized stopping times for which $\pi$ is the distribution of the
stopping state.
In symbols,
%
%
\begin{equation}
\label{LWmixdef} T_{\mathrm{mix}}(X):= \inf\EE S,
\end{equation}
where the infimum is taken over randomized stopping times $S$ such that
the distribution of $X_S$ is $\pi$. For computing $T_{\mathrm{mix}}(X)$, a
very useful theorem from~\cite{LovaszWinkler2} asserts that a
randomized stopping time $S$ achieves the minimum in (\ref{LWmixdef})
if and only if it has a halting state, that is, a state $x$ such that
if $X_t = x$ then (almost surely) $S \leq t$. We will use this result
to compute $T_{\mathrm{mix}}(X)$ for any irreducible birth-and-death chain
in Theorem~\ref{TTmixBD}, but first we state a lemma about expected hitting
times for birth-and-death chains.
%
%
\begin{lemma}
\label{LKTET}
For an irreducible birth-and-death chain on $\Xc= \{0,\ldots, n\}$ (in
discrete or continuous time) with stationary
distribution $\pi$ and initial state $0$, let $T$ denote the hitting
time of state $n$.\vspace*{8pt}

\textup{(a)} In discrete time, denote the birth probability from state $i$
by $p_i$. Then
\[
\EE T = \sum_{i = 0}^{n - 1} \frac{1}{\pi_i p_i}
\sum_{k = 0}^i \pi_k.
\]

\textup{(b)} In continuous time, denote the birth rate from state $i$ by
$\gl_i$. Then
\[
\EE T = \sum_{i = 0}^{n - 1} \frac{1}{\pi_i \gl_i}
\sum_{k =
0}^i \pi_k.
\]
\end{lemma}
\begin{pf}
Each assertion is easily established, and each follows immediately from
the other; for (b),
see, for example,~\cite{MR0356197}, Chapter 4, Problem 22.
\end{pf}
%
%
\begin{theorem}
\label{TTmixBD}
Let $X$ be an irreducible (discrete-time) birth-and-death chain on $\Xc
= \{0,\ldots, n\}$ with stationary p.m.f. (resp., c.d.f.) $\pi$ (resp., $H$)
and initial state~$0$. Then
\[
T_{\mathrm{mix}}(X) = \sum_{i = 0}^{n - 1}
\frac{H_i (1 - H_i)}{\pi_i p_i}.
\]
\end{theorem}
\begin{pf}
Let us use the \textit{naive rule} $S$ as our randomized stopping time:
choose $j$ randomly according to $\pi$, and then let $S$ be the hitting
time of $j$. Obviously the stopping distribution is $\pi$, as required.
Moreover, the
state $j$ must be hit en route to $n$; hence $n$ is a halting state
and $S$ achieves the minimum
at (\ref{LWmixdef}).

To compute $T_{\mathrm{mix}}(X) = \EE S$, we first note that Lemma \ref
{LKTET}(a) yields (easily) corresponding formulas for the expected
value of the hitting time $T_j$ of each state~$j$:
\[
\EE T_j = \sum_{i = 0}^{j - 1}
\frac{H_i}{\pi_i p_i}.
\]
Therefore
\[
T_{\mathrm{mix}}(X) = \sum_{j = 0}^n
\pi_j \EE T_j = \sum_{j = 0}^n
\pi_j \sum_{i = 0}^{j - 1}
\frac{H_i}{\pi_i p_i} = \sum_{i = 0}^{n - 1}
\frac{H_i (1 - H_i)}{\pi_i p_i}
\]
as desired.
\end{pf}
%
%
\begin{Remark}
\label{RTmixvssep} (a) The Lov\'asz--Winkler theory of mixing times and
the statement and proof of Theorem~\ref{TTmixBD} all carry over
routinely to the ``continuized'' chain which evolves in the same\vadjust{\goodbreak} way as
the given discrete-time chain but with independent exponential random
times with mean $1$ replacing unit times. In particular, the value of
$T_{\mathrm{mix}}(X)$ remains unchanged under continuization of an
irreducible discrete-time birth-and-death chain $X$ with initial state
$0$.

(b) By a theorem of Aldous and Diaconis
\cite{MR88d60175}, Proposition 3.2, in discrete time and a theorem of
Fill~\cite{MR1183166}, Theorem 1.1, in continuous time, any ergodic
finite-state Markov chain $X$ (regardless of initial distribution) has
a fastest (i.e., stochastically minimal) strong stationary time $T$
satisfying $\PP(T > t) = \operatorname{sep}(t)$ for every $t$
(restricted to integer
values for a discrete-time chain). If the state space is partially
ordered with bottom element $\zh$ and top element $\oneh$ and the
chain $X$ starts in $\zh$, and if the time-reversed kernel $\Kt$ is
monotone, then $\oneh$ is a halting state for any such $T$; to see
this, observe that
\begin{eqnarray*}
\PP(X_t = \oneh, T > t) &=& \PP(X_t = \oneh) - \PP(T \leq
t, X_t = \oneh)
\\
&=& \pi_{\oneh} \biggl[ \frac{K^t(\zeroh, \oneh)}{\pi_{\oneh}} - \bigl(1
-
\operatorname{sep}(t)\bigr) \biggr]
\\
&=& \pi_{\oneh} \biggl[ \min_i \frac{K^t(\zeroh, i)}{\pi_i} - \bigl(1 -
\operatorname{sep} (t)\bigr) \biggr] = 0,
\end{eqnarray*}
where $\pi$ is the stationary distribution and the penultimate equality
follows from the monotonicity of $\Kt^t$.

Now consider an ergodic birth-and-death chain $X$ (in discrete or
continuous time) on $\Xc= \{0,\ldots, n\}$ with stationary
distribution $\pi$ and initial state~$0$. In the discrete-time case,
assume that the chain is monotone; this is automatic in continuous time
by a simple and standard coupling argument. Then a fastest (i.e.,
stochastically minimal) strong stationary time $T$ exists, and $n$ is a
halting state for any such $T$. It follows that $T_{\mathrm{mix}}(X) = \EE
T$ and thus Theorem~\ref{TTmixBD} also gives an expression for
$\EE T$, which equals
\[
\sum_{t = 0}^{\infty} \PP(T > t) = \sum
_{t = 0}^{\infty} \operatorname{sep}(t)
\]
in discrete time and equals
\[
\int_0^{\infty} \PP(T > t) \,dt = \int
_0^{\infty} \operatorname{sep}(t) \,dt
\]
in continuous time. This remark gives added import to the value of
$T_{\mathrm{mix}}(X)$ for any irreducible discrete-time birth-and-death
chain $X$ (whether monotone or not) with initial state $0$: it equals
the integral of separation for the continuized chain.

(c) Given a collection $\mathcal{C}$ of irreducible discrete-time
birth-and-death chains $Y$ with initial state $0$, suppose that $X \in
\mathcal{C}$ satisfies $X = \arg\min_{Y \in\mathcal{C}}
T_{\mathrm{mix}}(Y)$. In light of remark (b), one might wonder whether
the continuized chain corresponding to $X$ minimizes $\operatorname
{sep}(t)$ at \textit{every} time $t$ over all continuizations of
chains $Y \in\mathcal {C}$. Theorem~\ref{TpathLW}(b) provides a
counterexample. Indeed, it can be shown that if we compare the chain of
the form (\ref{pchoice}) but with $\theta_n$ changed to $(n - 1) / (2
n)$ with any other birth-and-death chain having initial state $0$ and
symmetric kernel $K$, then there exists $t_0 = t_0(K)$ such that
continuized separation at time $t$ is strictly smaller for the former
chain than for the latter for all $0 < t \leq t_0$.\footnote{Indeed, if
$Y$ and $Z$ are the discrete-time and continuized chain corresponding
to $K$, then, with $\pi$ denoting the uniform p.m.f., as $t \to0$ we find
\[
1 - \operatorname{sep}_Z(t) = \frac{\PP(Z_t = n)}{\pi_n} = e^{-t}
\frac{t^n}{n!} \frac
{\PP(Y_n = n)}{\pi_n} + o\bigl(t^{n + 1}\bigr) =
\frac{t^n}{n!} (n + 1) p_0 p_1 \cdots
p_{n - 1} + o\bigl(t^{n + 1}\bigr),
\]
and $p_0 p_1 \cdots p_{n - 1}$ is uniquely maximized subject to $p_{k
-1} + p_k \leq1$ for $k = 0,\ldots, n - 1$ by choosing $p_k = (n + 1)
/ (2 n)$ if $k$ is even and $p_k = (n - 1) / (2 n)$ if $k$ is odd.}
Likewise, in the ``ladder game'' discussed in Section~\ref{Sladder}
it is the
uniform chain, not the chain discussed there, that is ``best in
separation for small $t$'' in similar fashion.
\end{Remark}

We are now in position to determine, for given $\pi$, the
birth-and-death chain $X$ that minimizes
$T_{\mathrm{mix}}(X)$ among those having initial state $0$, stationary
distribution $\pi$ and no holding probability except at the endpoints
of the state space. Unlike in Section~\ref{SBD}, we do not need to restrict
to monotone kernels; and rather than assuming that $\pi$ is
log-concave, we assume instead that $\pi$ is nondecreasing. For the
case that $\pi$ is uniform, we will give later an argument that removes
the restriction about holding probabilities. [There are examples, such
as $\pi= \frac{1}{15}(1, 2, 4, 4, 4)$, showing that the restriction
cannot be removed in general.]
%
%
\begin{theorem}
\label{TFMMCBDLW}
Let $\Xc= \{0,\ldots, n\}$. Among all irreducible birth-and-death
chains $X$ having a given positive nondecreasing stationary p.m.f. $\pi$,
initial state $0$ and no holding probability except at $0$ and $n$,
there is a unique chain $X_{\pi}$ minimizing $T_{\mathrm{mix}}(X)$.
Moreover:\vspace*{8pt}

\textup{(a)} Let $a_i:= \sum_{j = 1}^i (-1)^{i - j} \pi_j$ for $i = 0,\ldots
, n - 1$. Define
\[
f(w):= \sum_{i = 0}^{n - 1} \frac{H_i (1 - H_i)}{(-1)^i w + a_i}.
\]
Then there exists a unique $w_{\pi}$ minimizing $f(w)$ over $w \in[0,
\pi_0]$, and\break $T_{\mathrm{mix}}(X_{\pi}) = f(w_{\pi})$.

\textup{(b)} The optimal chain $X_{\pi}$ has transition probabilities
\[
q_i = \frac{a_{i - 1} + (-1)^{i - 1} w_{\pi}}{\pi_i},\qquad r_i = 0,\qquad p_i =
\frac{a_i + (-1)^i w_{\pi}}{\pi_i} \qquad(i = 0,\ldots, n)
\]
with the exceptions $q_0 = 0$, $r_0 = 1 - p_0$, $r_n = 1 - q_n$ and
$p_n = 0$.
\end{theorem}
\begin{pf}
We begin by noting
that birth-and-death kernels with stationary distribution $\pi$ (in
complete generality, irrespective of holding probabilities or
nondecreasing\-ness of $\pi$) are in one-to-one correspondence with
nonnegative sequences
$\ww= (w_{-1}, w_0,\ldots, w_n)$ satisfying $w_{-1} = 0 = w_n$ and
%
%
\begin{equation}
\label{wineq} w_{i - 1} + w_i \leq\pi_i \qquad(i
= 0,\ldots, n),
\end{equation}
the correspondence being $w_i = \pi_i p_i = \pi_{i + 1} q_{i + 1}$, $i
= 0,\ldots, n - 1$. The proof is easy, and the correspondence gives
\[
r_i = 1 - q_i - p_i = 1 - \frac{w_{i - 1} + w_i}{\pi_i}
\qquad(i = 0,\ldots, n)
\]
for the holding probabilities. In this $\ww$-parameterization,
Theorem~\ref{TTmixBD} gives
%
%
\begin{equation}
\label{TmixBD} T_{\mathrm{mix}} = \sum_{i = 0}^{n - 1}
\frac{H_i (1 - H_i)}{w_i}.
\end{equation}

The constraint $r_i = 0$ for $i = 0,\ldots, n - 1$ is precisely the
constraint that equality holds in (\ref{wineq}) for
$i = 1,\ldots, n - 1$. Then we must have $w:= w_0 \in[0, \pi_0]$ and
\[
w_i = (-1)^i w + a_i \qquad(i = 0,\ldots, n
- 1).
\]
It follows from the assumption that $\pi$ is nondecreasing that these
$w_i$'s are indeed all nonnegative [and all positive if $w \in(0, \pi
_0)$]. This proves the theorem, because $f$ is continuous on $[0, \pi
_0]$ and both finite and strictly
convex\footnote{In the general setting of (\ref{TmixBD}), $T_{\mathrm{mix}}$
is a strictly convex function on a nonempty convex domain (an
intersection of half-spaces) of arguments $\ww$ and so has a unique
minimum. The optimal $\ww$ is on the boundary of the domain; more
specifically, for every $i = 0,\ldots, n - 1$, if the optimal $\ww$
does not lie on the hyperplane delimiting the $i$th half-space (\ref
{wineq}), then it lies on the $(i + 1)$st such hyperplane.}
on $(0, \pi_0)$.
\end{pf}

We now specialize to the case of uniform $\pi$, removing the
restriction on holding from Theorem~\ref{TFMMCBDLW} and solving explicitly
for the value $w$ in Theorem~\ref{TFMMCBDLW}(a). We find it somewhat
surprising that the chain minimizing $T_{\mathrm{mix}}$ is \textit{not} the
uniform chain whenever $n \geq3$ is odd.
%
%
\begin{theorem}
\label{TpathLW}
Consider the problem of minimizing $T_{\mathrm{mix}}$ among all
birth-and-death chains on $\Xc= \{0,\ldots, n\}$ with initial
state $0$ and symmetric kernel.\vspace*{8pt}

\textup{(a)} If $n \geq2$
is even, then the uniform chain is the unique minimizing chain.

\textup{(b)} If $n$ is odd, then
%
%
\begin{equation}
\label{pchoice} p_k = %
\cases{ 1 - \theta_n,
&\quad if $k$ is even
\cr
\theta_n, &\quad if $k$ is odd} %
\qquad(k = 0,\ldots, n - 1)\vadjust{\goodbreak}
\end{equation}
gives the unique minimizing chain, where for any $m$ we define
%
%
\begin{equation}
\label{theta} \theta_{m - 1}:= \tfrac{1}{6} \bigl[ \sqrt{
\bigl(m^2 + 2\bigr) \bigl(m^2 - 4\bigr)} -
\bigl(m^2 - 4\bigr) \bigr].
\end{equation}
\end{theorem}

We have written the formula for $\theta_{m - 1}$ rather than that for
$\theta_n$ because it is simpler to write.
%
%
\begin{Remark}
Although the uniform chain is not optimal when $n$ is odd, it is nearly
optimal, since $\theta_n$ has the asymptotics
\[
\theta_n = \tfrac{1}{2} - \tfrac{3}{4} n^{-2}
+ O\bigl(n^{-3}\bigr)\qquad\mbox{as $n \to\infty$}
\]
and the value of $T_{\mathrm{mix}}$ (recall Theorem~\ref{TTmixBD}) for $p_k
\equiv
1 / 2$ is $\frac{1}{3} n^2 + n + \frac{2}{3}$, only slightly larger
than the optimal value $\frac{1}{3} n^2 + n + \frac{2}{3} - \frac{3}{4}
n^{-2} + O(n^{-3})$.
\end{Remark}
\begin{pf*}{Proof of Theorem~\ref{TpathLW}}
Recall Theorem~\ref{TTmixBD}; thus the goal is to minimize
\[
f({\mathbf p}):= \sum_{k = 0}^{n - 1}
\frac{(k + 1)(n - k)}{p_k}
\]
over vectors ${\mathbf p}= (p_0,\ldots, p_{n - 1})$ that are nonnegative
(we will not repeat this nonnegativity condition below) and satisfy
%
%
\begin{equation}
\label{pinequ} p_{k - 1} + p_k \leq1\qquad\mbox{for $k = 0,\ldots, n$},
\end{equation}
where $p_{-1} = 0 = p_n$.
The objective function $f({\mathbf p})$ is strictly convex in ${\mathbf
p}$ (by strict
convexity of $x \mapsto x^{-1}$). Hence there is a unique minimizer,
and because $(p_{n - 1},\ldots, p_0)$ is clearly a minimizer if
$(p_0,\ldots, p_{n - 1})$ is, the unique minimizer is of the form
\[
(p_0,\ldots, p_{(n / 2) - 1}, p_{(n / 2) - 1},\ldots,
p_0),
\]
if $n$ is even and of the form
\[
(p_0,\ldots, p_{(n - 3) / 2}, p_{(n - 1) / 2}, p_{(n - 3) / 2},\ldots,
p_0),
\]
if $n$ is odd. We now break into the two cases.\vspace*{8pt}

(a) For $n$ even, we seek equivalently to minimize
\[
f({\mathbf p}) = 2 \sum_{k = 0}^{(n / 2) - 1}
\frac{(k + 1) (n - k)}{p_k}
\]
subject to
\[
p_{k - 1} + p_k \leq1\qquad\mbox{for $k = 0,\ldots, (n / 2)$}.
\]
(Note that the last of these conditions is $p_{(n / 2) - 1} \leq1 /
2$.)\vadjust{\goodbreak}

We claim (by induction on $K$) for $1 \leq K \leq(n / 2) - 1$ that the
minimizer of $\sum_{k = 0}^K \frac{(k + 1) (n - k)}{p_k}$ subject to
(nonnegativity and)
$p_{k - 1} + p_k \leq1$ for $k = 0,\ldots, K$ and $p_K \leq1 / 2$ is
$p_k \equiv1 / 2$.

For the basis $K = 1$ of the induction, we seek to minimize
\[
\frac{n}{p_0} + \frac{2 (n - 1)}{p_1}
\]
subject to $p_0 + p_1 \leq1$ and $p_1 \leq1 / 2$. Clearly we should
take $p_0 = 1 - p_1$ (regardless of $p_1$), and then we need to minimize
\[
\frac{n}{1 - p_1} + \frac{2 (n - 1)}{p_1}
\]
subject to $p_1 \leq1 / 2$. Because $2 (n - 1) \geq n$ (i.e., $n \geq
2$), the minimizer is $p_1 = 1 / 2$ (and then $p_0 = 1 / 2$).

We now proceed to the induction step to move from $K - 1$ to $K$. To
minimize, clearly we should take
$p_K = \min\{ 1 / 2, 1 - p_{K - 1}\}$. The remainder of the proof
for $n$ even then breaks into two cases.

\textsc{Case} 1. If $p_{K - 1} \geq1 / 2$, then we take $p_K = 1 - p_{K -
1}$ and our goal is to minimize
\[
\sum_{k = 0}^{K - 2} \frac{(k + 1) (n - k)}{p_k} +
\frac{K (n - (K -
1))}{p_{K - 1}} + \frac{(K + 1) (n - K)}{1 - p_{K - 1}}
\]
subject to $p_{k - 1} + p_k \leq1$ for $0 \leq k \leq K - 1$ and
(because this is case 1) $p_{K - 1} \geq1 / 2$. Because
$(K + 1) (n - K) \geq K (n - (K - 1))$ and we have the restriction
$p_{K - 1} \geq1 / 2$, we should set $p_{K - 1}$ as small as possible,
namely, $p_{K - 1} = 1 / 2$, and then we seek to minimize
\[
\sum_{k = 0}^{K - 2} \frac{(k + 1) (n - k)}{p_k} +
\frac{K (n - (K -
1))}{p_{K - 1}} + \frac{(K + 1) (n - K)}{1 / 2}
\]
subject to $p_{k - 1} + p_k \leq1$ for $0 \leq k \leq K - 1$ and $p_{K
- 1} = 1 / 2$. Clearly the minimum value here is at least as large as
the minimum value if we relax the last constraint to $p_{K - 1} \leq1
/ 2$. But then by induction the minimum value is achieved by setting
$p_k \equiv1 / 2$. This completes the proof in case 1.

\textsc{Case} 2. If $p_{K - 1} \leq1 / 2$, then we set $p_K = 1 / 2$ and
the goal is to minimize
\[
\sum_{k = 0}^{K - 1} \frac{(k + 1) (n - k)}{p_k} +
\frac{(K + 1) (n -
K)}{1 / 2}
\]
subject to $p_{k - 1} + p_k \leq1$ for $0 \leq k \leq K$ and $p_{K -
1} \leq1 / 2$. But then again by induction the minimum value is
achieved by setting $p_k \equiv1 / 2$. This completes the proof in
case 2, and thereby completes the proof of part (a).

(b) For $n$ odd, suppose without loss of generality that $n \geq3$. We
first prove that the optimum is again attained for a chain that
satisfies equality in
condition (\ref{pinequ}) at interior points $k$ of the state space:
%
%
\begin{equation}
\label{pkeq} p_{k - 1} + p_k = 1\qquad\mbox{for $k = 1,\ldots, n
- 1$}.
\end{equation}
Recall that the minimizing ${\mathbf p}$ is unique and symmetric. Hence,
considering the holding probability
$r_k:= 1 - p_{k - 1} - p_k$ at state $k$, it suffices to show that
there is an optimizing chain with $r_k = 0$ for
$1 \leq k \leq(n - 1) / 2$.

We proceed by contradiction. We show that there exists ${\mathbf p}'$
satisfying (\ref{pinequ}) and $f({\mathbf p}') < f({\mathbf p})$ in
each of the
following three cases which, allowing arbitrary $k \in\{1,\ldots, (n -
1) / 2\}$, exhaust all possibilities where $r_k > 0$ for some $1 \leq k
\leq(n - 1) / 2$:
\begin{longlist}
\item$r_k > 0$ and $r_{k - 1} > 0$;
\item$r_k > 0$ and $r_{k - 1} = 0$ and $p_k \geq1 / 2$;
\item$p_k < 1/2$, and $k$ is the largest value $j$ in $\{1,\ldots, (n
- 1) / 2\}$ such that $r_j > 0$.
\end{longlist}

In case (i), let
\[
p'_{k - 1}:= p_{k - 1} + \min\{r_{k - 1},
r_k\}
\]
and $p'_j:= p_j$ otherwise.

In case (ii),
first note that $k \geq2$; indeed, were we to have $k = 1$, then (by
our assumption) $r_0 = 0$ and so $p_0 = 1$; but then $p_1 = 0$, and
such a ${\mathbf p}$ clearly does not minimize $f({\mathbf p})$.
Next, because $p_k \geq1 / 2$ we must have $p_{k - 1} < 1 / 2$
(because $r_k > 0$) and thus
$p_{k-2} > 1 / 2$ (because $r_{k - 1} = 0$).
We can then let
\[
p'_{k-1}:= p_{k-1} + \varepsilon,\qquad
p'_{k-2}:= p_{k-2} - \varepsilon
\]
for suitably small $\varepsilon> 0$, and $p'_j:= p_j$ otherwise.
Since $k \leq(n - 1) / 2$, we know $k (n + 1- k) > (k - 1) (n + 2 -
k)$, so the derivative of $f({\mathbf p})$ in the direction of the vector
$\delta_{k - 1} - \delta_{k - 2}$ is negative and $f({\mathbf p}') <
f({\mathbf p})$.

In case (iii) we have $p_{k + 2 i} = p_k$ for $0 \leq i \leq\frac{n -
1}{2} - k$,
and $p_{k + 2 i - 1} = 1 - p_k$ for $1 \leq i \leq\frac{n - 1}{2} -
k$. We form ${\mathbf p}'$ by changing these values to
$p'_{k + 2 i}:= p_k + \varepsilon$ and $p'_{k + 2 i - 1}:= 1 - p_k -
\varepsilon$ for suitably small $\varepsilon> 0$ and setting $p'_j:= p_j$
otherwise. We see that $f({\mathbf p}') < f({\mathbf p})$ if the
derivative with
respect to $p_k$ of the following expression is negative for all $p_k <
1 / 2$:
\[
\frac{1}{p_k} \sum_{i = 0}^{({n - 1})/{2} - k} (k + 2 i
+ 1) (n - k - 2 i) + \frac{1}{1 - p_k}\sum_{i = 1}^{({n - 1})/{2} - k}
(k + 2i ) (n + 1 - k - 2 i);
\]
and that is true if (and only if) the first sum is at least as large as
the second. Indeed, the first sum \textit{is} larger than the second:
\begin{eqnarray*}
&&
\sum_{i = 0}^{({n - 1})/{2} - k} (k + 2 i + 1) (n - k - 2 i)
- \sum_{i = 1}^{({n - 1})/{2} - k} (k + 2i ) (n + 1 - k - 2
i)
\\
&&\qquad= (k + 1) (n - k) + \sum_{i = 1}^{({n - 1})/{2} - k} (n -
2 k - 4 i)
\\
&&\qquad= k (n - k) + \frac{1}{2} (n + 1) > 0.
\end{eqnarray*}

Since we have established constraint (\ref{pkeq}), \textit{every}
feasible vector ${\mathbf p}$ is of the form
\[
p_k \equiv%
\cases{ 1 - \theta, &\quad if $k$ is even,
\cr
\theta, &\quad if $k$ is odd, } %
\]
so we need only verify that the choice $\theta= \theta_n$ as defined
at (\ref{theta}) is optimal. Indeed, writing $r = (n - 1) / 2$ we have
\begin{eqnarray*}
a_n:= \sum_{0 \leq j \leq r} (2 j + 1) (n - 2 j) &=&
\frac{1}{12} (n + 1) \bigl(n^2 + 2n + 3\bigr),
\\
b_n:= \sum_{1 \leq j \leq r} (2 j) (n - 2j + 1) &=&
\frac{1}{12} (n + 1) (n - 1) (n + 3),
\end{eqnarray*}
and then the optimal choice of $\theta$, minimizing $\frac{a_n}{1 -
\theta} + \frac{b_n}{\theta}$, is $\theta_n$ given by
\[
\theta_n = ( 1 + \sqrt{a_n / b_n}
)^{-1}.
\]
After a little bit of computation, we find that $\theta_n$ is given in
accordance with equation (\ref{theta}).
\end{pf*}

\section{\texorpdfstring{A ``ladder'' game.}{A ``ladder'' game}}\label{Sladder}
In this section we discuss a simple ``ladder'' game, where the class of
kernels considered is a certain subclass of the symmetric
birth-and-death kernels considered in Section~\ref{Sunif}. Our treatment
involves finding the kernel that minimizes
the Lov\'{a}sz--Winkler mixing time $T_{\mathrm{mix}}$.
This particular kernel is \textit{not} one that had previously been
considered as a candidate for ``fastest.''

Lange and Miller~\cite{LM} discusses a ``ladder'' game and several
contexts, including an old Japanese scheme for choosing a spouse's
Christmas gift from a list of desired items, in which it arises. We
refer the reader to~\cite{LM} for details. A class of Markov chains
that arise in modeling the ladder game (see ``Model One'' in
\cite{LM}, Section 5) have the permutation group on $\{0,\ldots, n\}$ as
state space and moves that transpose items in adjacent positions; write
$p_i$ for the probability that the positions chosen are $i$ and $i +
1$, so that
%
%
\begin{equation}
\label{ladder} p_0 + p_1 + \cdots+ p_{n - 1} =
1.
\end{equation}
We will refer to (\ref{ladder}) as the ``ladder condition.'' If we
follow the movement of only a single item (this is ``Model Two: The
path of a single marcher as a random walk among the columns of the
ladder'' in~\cite{LM}, Section 7, esp. Figure 9), then we have
precisely the class of symmetric birth-and-death kernels considered in
our path-problem of Section~\ref{Sunif}, but now subject to the
ladder condition.
From~\cite{LM} (Section~8: How many rungs is enough?) we have the following
quote (with notation adjusted slightly to match that of
Section~\ref{Sunif}):

\begin{quote}
We suspect (but have not shown) that for any $n$, the rate of
convergence is maximized when rung placement is uniform. That is, the
absolute value of the largest small eigenvalue is minimized when $p_i =
1 / n$ for $i = 0, 1,\ldots, n - 1$.
\end{quote}

\noindent (Here ``largest small eigenvalue'' means the eigenvalue of the kernel
with largest absolute value strictly less than $1$---what is
called ``SLEM''
in~\cite{MR2124681,MR2202924,MR2515797}.) The authors of~\cite{LM}
base their suspicion on calculations for
$n = 2$, for which their conjecture is indeed true.

The corresponding continuous-time problem has been studied by
Fied\-ler~\cite{MR1034418} and, in a somewhat more general setting, by
Sun et al. in~\cite{MR2278445}, Example 5.2. The result is that, among
all continuous-time symmetric birth-and-death chains on $\{0,\ldots,
n\}
$, started from $0$, with birth rates $p_i$ satisfying the ladder
condition~(\ref{ladder}), the one which is fastest-mixing in the sense
of minimizing relaxation time has $p_i$ proportional to $(i + 1) (n -
i)$. It can be shown
that these weights also\vspace*{1pt} uniquely minimize SLEM in discrete time, so the
conjecture in~\cite{LM} is false for
every $n \geq3$.\footnote{At the end of their Section 8, the authors
of~\cite{LM} also wonder, based on results for $n = 2$, whether it
might be the case for all $n$ that, except for multiplicities, the
eigenvalues are the same for the permutation chain as for the
single-marcher chain. This is seen to be false by the discussion in
\cite{MR2629990}, Section 1.4. But the main theorem of~\cite{MR2629990}
does establish that the second-largest eigenvalues of the two chains agree.}

One might now suspect that these parabolic weights provide a FMMC
(subject to the ladder condition) in a variety of senses, at least for
chains (as henceforth assumed) starting in state $0$. However, working
in discrete time, it
is clear (a) from reviewing the discussion in
Section~\ref{SSmono} that there is no bottom element with respect
to $\preceq
$ for monotone chains satisfying the ladder condition and
(b) from Remark~\ref{R2stepCI} that there is no bottom element
in $\preceq$
for squares of ladder-condition birth-and-death kernels. Further, it
can be shown, switching
to continuous time to match the setting of~\cite{MR2278445} and in
order to bring standard techniques to bear (it is well known that all
birth and death chains in continuous time are monotone), that there is
no ladder-condition birth-and-death chain minimizing separation at
every time. Theorem~\ref{Tladder} implies that the integral of separation
over all times is minimized by weights $p_i$ proportional to the square
roots $\sqrt{(i + 1) (n - i)}$ of the weights minimizing SLEM.
%
%
\begin{theorem}
\label{Tladder}
For each discrete-time symmetric birth-and-death chain with state space
$\{0,\ldots, n\}$, initial state $0$ and birth probabilities\vadjust{\goodbreak} ${\mathbf
p} =(p_i)$ satisfying the ladder condition (\ref{ladder}), let $f({\mathbf p})$
denote its Lov\'{a}sz--Winkler mixing time $T_{\mathrm{mix}}$. Then
the uniquely optimal (i.e., minimizing) choice of ${\mathbf p}$ is to take
$p_i$ proportional to $\sqrt{(i + 1) (n - i)}$.
\end{theorem}

Theorem~\ref{Tladder} is an immediate consequence of the following corollary
to the proof of Theorem~\ref{TFMMCBDLW}, taking $\pi$ to be uniform and $c$
to be
$1 / n$.
%
%
\begin{corollary}
\label{Cmin}
Over all discrete-time birth-and-death chains on $\{0,\break\ldots, n\}$
(started at $0$) with given stationary
distribution $\pi$ (having c.d.f. $H$) and
\[
\sum_{k = 0}^{n - 1} \pi_k
p_k = c \in\Bigl(0, \min_i \pi_i\Bigr],
\]
the mixing time $T_{\mathrm{mix}}$ of the chain is minimized by the choice
\begin{eqnarray*}
q_k &\equiv&\frac{c \sqrt{H_{k - 1} (1 - H_{k - 1})}}{\pi_k \sum_j
\sqrt{H_j (1 - H_j)}},\qquad p_k \equiv
\frac{c \sqrt{H_k (1 - H_k)}}{\pi_k \sum_j \sqrt{H_j
(1 - H_j)}},\\
r_k &\equiv&1 - q_k -
p_k,
\end{eqnarray*}
and the minimized value is
\[
T_{\mathrm{mix}} = c^{-1} \Biggl[ \sum_{k = 0}^{n - 1}
\sqrt{H_k (1 - H_k)} \Biggr]^2.
\]
\end{corollary}
\begin{pf}
As demonstrated in the proof of Theorem~\ref{TFMMCBDLW}, the goal is
to minimize
\[
T_{\mathrm{mix}} = \sum_{i = 0}^{n - 1}
\frac{H_i (1 - H_i)}{w_i}
\]
over nonnegative sequences $(w_{-1}, w_0,\ldots, w_n)$ satisfying
$w_{-1} = 0 = w_n$ and
%
%
\begin{equation}
\label{ignore} w_{i - 1} + w_i \leq\pi_i \qquad(i
= 0,\ldots, n)
\end{equation}
and $\sum_{k = 0}^{n - 1} w_k = c$. Ignoring the constraint (\ref
{ignore}), the optimal choice of
the weights $w_i$ is clear, namely, $w_i \equiv\pi_i p_i$ with $p_i$
as asserted in the statement of the theorem. But then (\ref{ignore}) is
automatically satisfied because we assume $c \in(0, \min_i \pi_i]$.
Evaluation of the objective function at the optimizing kernel gives the
optimized value of $T_{\mathrm{mix}}$.
\end{pf}
%
%
\begin{Remark}
\label{laddercomparison}
Let
$n \to\infty$. For the optimal kernel of Theorem~\ref{Tladder} we have
$T_{\mathrm{mix}} \sim\frac{\pi^2}{64} n^3$, whereas for both $p_i \equiv1
/ n$ (the guess for optimality in~\cite{LM}) and the choice $p_i
\propto(i + 1) (n - i)$ minimizing SLEM we have $T_{\mathrm{mix}} = \frac
{1}{6} n (n + 1) (n + 2) \sim\frac{1}{6} n^3$.
\end{Remark}

\section{\texorpdfstring{Can extra updates delay mixing? (No, subject to
positive correlations.)}{Can extra updates delay mixing? (No, subject to
positive correlations)}}\label{SPW}

Can extra updates delay mixing? This question is the title of a
paper~\cite{PW2011} by Yuval Peres and Peter Winkler; see also
Holroyd~\cite{holroyddelay} for counterexamples.\vadjust{\goodbreak} Peres and Winkler
show that the answer is no, for total variation distance, in the
setting of monotone spin systems, generalized by replacing the set of
spins $\{0, 1\}$ by any linearly ordered set. (We review relevant
terminology below.) In Theorem~\ref{TPW} we recapture and extend
their result
using comparison inequalities by showing that $K_v \preceq I$ for any
kernel $K_v$ that updates a single site $v$, that is, that the identity
kernel [as for the monotone birth-and-death example, see Remark \ref
{RBDmono}(a)] only slows mixing (when the initial p.m.f. has
nonincreasing ratio with respect to the stationary p.m.f.)---because
then, noting reversibility and stochastic monotonicity of each $K_v$
and applying Proposition~\ref{PCI}, for any $v_1,\ldots, v_t$ the product
$K_{v_1} \cdots K_{v_t}$ increases in $\preceq$ by deletion of any
$K_{v_i}$. The comparison inequality $K_v \preceq I$ holds in the more
general setting of a partially ordered set of ``spins,'' subject to the
following restriction: starting with distribution $\pi$ and a site $v$
and conditioning on the spins at all sites other than $v$, the
conditional law of the spin at $v$ should have positive correlations
(as, of course,
does any distribution on a \textit{linearly} ordered set).

\subsection{\texorpdfstring{Positive correlations.}{Positive correlations}}\label{SSposcorr}
Recall that a p.m.f. $\pi$ on a finite partially ordered set $\Xc$ is
said to
\textit{have positive correlations} if (in the notation of Section~\ref{SCI})
\[
\langle f, g \rangle\geq\langle f, 1 \rangle\langle g, 1 \rangle
\]
for every $f, g \in\MMc$, and that if $S$ is \textit{linearly} ordered
then (by ``Chebyshev's other inequality;'' see, e.g.,
\cite{PeresWinkler}, Lemma 16.2) \textit{all} probability measures have
positive correlations. The connection with comparison inequalities is
the following simple lemma, in relation to which we note that both
$K_{\pi}$ and $I$ are stochastically monotone kernels possessing
stationary distribution $\pi$.
%
%
\begin{lemma}
\label{Lposcorr}
A p.m.f. $\pi$ on a finite partially ordered set $\Xc$ has positive
correlations if and only if $K_{\pi} \preceq I$, where
$K_{\pi}$ is the trivial kernel that jumps in one step to $\pi$ and $I$
is the identity kernel.
\end{lemma}
\begin{pf}
Since for any $f$ and $g$ we have
\[
\langle K_{\pi} f, g \rangle= \bigl\langle\langle f, 1 \rangle, g \bigr
\rangle= \langle f, 1 \rangle\langle g, 1 \rangle
\]
and $\langle I f, g \rangle= \langle f, g \rangle$, the lemma is proved.
\end{pf}
%
%
\begin{proposition}
\label{Pposcorr}
Let $\pi$ be a p.m.f. on a finite partially ordered set. Partition~$\Xc$,
suppose that a given
kernel $K$ on $\Xc$ is a direct sum
[as in Proposition~\ref{PCIbasic}\textup{(c)}] of trivial kernels $K_i$ (as in
Lemma~\ref{Lposcorr}) on the cells of the partition, and suppose
that $\pi$ conditioned to each cell has positive correlations. Then
$K \preceq I$.
\end{proposition}
\begin{pf}
Simply combine Lemma~\ref{Lposcorr} and Proposition~\ref{PCIbasic}(c).
\end{pf}

\subsection{\texorpdfstring{Monotone spin systems.}{Monotone spin systems}}\label{SSspin}

Our setting is the following. We are given a finite graph $G = (V, E)$
and a finite partially ordered set $S$ of ``spin values.'' A \textit{spin
configuration} is an assignment of spins to vertices (sites), and our
state space is the set $\Xc$ of all configurations. We are given a
p.m.f. $\pi$ on $\Xc$ that is \textit{monotone} in the sense that, when we
start with $\pi$ and any site $v$ and condition on the spins at all
sites other than $v$, the conditional law of the spin
at $v$ is monotone in the conditioning spins. We recover
and (modestly) extend the Peres--Winkler result by means of the
following theorem, which (i) allows somewhat more general $S$ and
(ii) encompasses---by means of Proposition~\ref{PCIdom},
Corollary~\ref{CTVsepL2}(a) and (b),
and Remark~\ref{RL2}---separation and $L^2$-distance as well as TV.
%
%
\begin{theorem}
\label{TPW}
Fix a site $v$, and suppose that the conditional distributions
discussed in the preceding paragraph all have positive correlations.
Let $K_v$ be the (stochastically monotone) Markov kernel for update at
site $v$ according to the conditional distributions discussed. Then we
have the comparison inequality \mbox{$K_v \preceq I$}.
\end{theorem}
\begin{pf}
Say that two configurations are equivalent if they differ at most in
their spin at $v$, and let $[x]$ denote the equivalence class
containing a given configuration $x$. Then $K_v$ is given by
\[
K_v(x, y) = {\mathbf1}\bigl(y \in[x]\bigr) \frac{\pi(y)}{\pi([x])}.
\]
This $K_v$ is the direct sum of the trivial kernels (as in Lemma \ref
{Lposcorr}) on each equivalence class. Further, each class is
naturally isomorphic as a partially ordered set to $S$ and so has
positive correlations.
It is well known and easily checked that $K_v$ is stochastically
monotone, so
the theorem is an immediate consequence of Proposition~\ref{Pposcorr}.
\end{pf}
%
%
\begin{Remark}[(Random vs. systematic site updates)]
\label{randomvssystematic}
It follows [from Theorem~\ref{TPW} and Proposition \ref
{PCIbasic}(b)] for monotone spin
systems with (say) linearly ordered $S$ that, when the chains start
from a common p.m.f. having nonincreasing ratio relative to $\pi$, the
``systematic site updates'' chain with kernel $K_{\mathrm{syst}}:= K_{v_1}
\cdots K_{v_{\nu}}$ (for any ordering $v_1,\ldots, v_{\nu}$ of the
sites $v \in V$) mixes faster in TV, sep, and $L^2$ than does the
``random site updates'' chain with kernel $K_{\mathrm{rand}}:= \sum_{v \in
V} p_v K_v$ [for any p.m.f. ${\mathbf p}= (p_v)_{v \in V}$ on $V$]. This is
because (recalling the paragraph preceding
Proposition~\ref{PCIbasic}) the reversible kernel $K_{\mathrm{rand}}$ is
stochastically monotone, as are $K_{\mathrm{syst}}$ and its time-reversal,
and $K_{\mathrm{syst}} \preceq K_{\mathrm{rand}}$.
(The explanation for the comparison here is that (as noted in the first
paragraph of this section)
$K_{\mathrm{syst}} \preceq K_v$ for each $v \in V$ and [by Proposition \ref
{PCIbasic}(b)] the relation $\preceq$ on $\Kc$ is preserved under
mixtures.) It is important to keep in mind here that one ``sweep'' of
the sites using $K_{\mathrm{syst}}$ is counted as only one Markov-chain step.\vadjust{\goodbreak}

There is a very weak ordering in the opposite direction:
$K_{\mathrm{rand}}^{\nu} \preceq p K_{\mathrm{syst}} + (1 - p) I$, with $p:=
\prod_{v \in V} p_v$.
\end{Remark}

\subsection{\texorpdfstring{Extra updates do not delay mixing: Card-shuffling.}{Extra updates do not delay mixing: Card-shuffling}}\label{SScards}
The following
card-shuffling Markov chain, which has been studied quite a bit
(see~\cite{MR2135733} and references therein)
in the time-homogeneous ``random updates'' case where update positions
are chosen independently and uniformly,
is another example where comparison inequalities can be used to show
that extra updates do not delay mixing.

Our state space is the set $\Xc$ of all permutations of $\{1,\ldots,
n\} $, and there is a parameter $p \in(0, 1)$. Given $i \in\{1,\ldots,
n - 1\}$, we can \textit{update} adjacent positions $i$ and $i + 1$ by
sorting (i.e., putting into natural order) the two cards (numbers) in
those positions with probability $p$ and ``anti-sorting'' them with the
remaining probability. Call the update kernel $K_i$. It is
straightforward to check that each $K_i$ is (i) reversible with respect
to $\pi$, where $\operatorname{inv}(x)$ is the number of inversions in
the permutation $x$ and $\pi(x)$ is proportional to $[(1 - p) /
p]^{\mbox{\scriptsize inv}(x)}$ [indeed, $K_i(x, \cdot)$ is the law of
a permutation drawn from $\pi$ but conditioned to agree with $x$ at all
positions other than $i$ and $i + 1$], and (ii) stochastically monotone
with respect to the Bruhat order on $\Xc$ (defined so that $x \leq y$
if $y$ can be obtained from $x$ by a sequence of anti-sorts of
\textit{not necessarily adjacent} cards).\footnote{To establish the
monotonicity of $K_i$, it is sufficient to consider initial states $x$
and $y$ where $y$ is obtained from $x$ by a single anti-sort of two not
necessarily adjacent cards and couple transitions from these states so
that the corresponding terminal states, call them $X_1$ and $Y_1$,
satisfy $X_1 \leq Y_1$. A~coupling that one can check works (by
considering various cases) is to make the \textit{same} decision, for
$x$ and for $y$, to sort or to anti-sort the cards in positions $i$ and
$i + 1$.}
%
%
\begin{theorem}
\label{Tcards}
Fix a position $i \in\{1,\ldots, n - 1\}$, and let $K_i$ be the Markov
kernel for update of positions $i$ and $i + 1$ as discussed in the
preceding paragraph. Then we have the comparison inequality
$K_i \preceq I$.
\end{theorem}
The proof of Theorem~\ref{Tcards} is essentially the same as
for Theorem~\ref{TPW}
and therefore is omitted.
The key is that the relevant equivalence classes now consist of only
two permutations each and are linearly ordered, therefore
having positive correlations.

\subsection{\texorpdfstring{A final example.}{A final example}}\label{SSfinal}

In a specific setting (linearly ordered state space and uniform
stationary distribution) we have $K \preceq I$ quite generally:
%
%
\begin{theorem}
\label{TtopI}
Let $\Xc$ be a linearly ordered state space. If $K$ is doubly
stochastic, then $K \preceq I$
(with respect to uniform $\pi$).
\end{theorem}
%
%
\begin{Remark}
(a) When $\pi$ is uniform, to say that a kernel $K$ is doubly
stochastic is precisely to say that $\pi$ is stationary
for $K$. If $K$ is symmetric, then Theorem~\ref{TtopI} applies. Thus
inserting a monotone symmetric kernel (or, more generally, a monotone
doubly stochastic kernel whose transpose is also monotone) in a list of
such kernels to be applied never slows mixing (by Proposition \ref
{PCI}, or the
more general Corollary~\ref{Cc1}, and the results of Section~\ref
{Smaj}) when the
initial p.m.f. is nonincreasing.

(b) If ``linearly ordered'' is relaxed to ``partially ordered'' in
Theorem~\ref{TtopI}, the result is not generally true, even for
monotone $K$. This
follows from Lemma~\ref{Lposcorr}, since there are
partially ordered sets for which the uniform distribution does not have
positive correlations.
\end{Remark}
\begin{pf*}{Proof of Theorem~\ref{TtopI}}
We must show that $\langle K f, g \rangle\leq\langle f, g \rangle$
when $f$ and $g$ are nonnegative and
belong to $\MMc$ (i.e., are nonincreasing) and (without loss of
generality) $f$ sums to $1$. It is a fundamental result in the theory
of majorization~\cite{MR81b00002} that $f$ majorizes $K f$ if $K$ is
doubly stochastic. Since $\Xc$ is linearly ordered and $f$ belongs
to $\MMc$, it follows that, regarded as p.m.f.'s, $f$ and $K f$ satisfy
$K f \geq f$ stochastically. Therefore, for $g \in\MMc$ we have
$\langle K f, g \rangle\leq\langle f, g \rangle$, as desired.
\end{pf*}

\section*{\texorpdfstring{Acknowledgments.}{Acknowledgments}}

We thank Persi Diaconis, Roger Horn, Chi-Kwong
Li, Guy Louchard and Mario Ullrich for helpful discussions; an
anonymous referee for helpful comments; and an anonymous Associate
Editor for alerting us to the reference~\cite{PW2011} (we had
previously cited their work as unpublished).



\printaddresses


\begin{thebibliography}{35}

\bibitem{MR88d60175}
\begin{barticle}[mr]
\bauthor{\bsnm{Aldous},~\bfnm{David}\binits{D.}} \AND
  \bauthor{\bsnm{Diaconis},~\bfnm{Persi}\binits{P.}}
(\byear{1987}).
\btitle{Strong uniform times and finite random walks}.
\bjournal{Adv. in Appl. Math.}
\bvolume{8}
\bpages{69--97}.
\bid{doi={10.1016/0196-8858(87)90006-6}, issn={0196-8858}, mr={0876954}}
\bptok{imsref}%
\end{barticle}
\endbibitem

\bibitem{AldousFill}
\begin{bmisc}[auto:STB|2012/12/05|11:57:16]
\bauthor{\bsnm{Aldous},~\bfnm{D.~J.}\binits{D.~J.}} \AND
  \bauthor{\bsnm{Fill},~\bfnm{James~Allen}\binits{J.~A.}}
\bhowpublished{Reversible Markov chains and random walks on graphs. Chapter
  drafts available at
  \texttt{\href{http://www.stat.Berkeley.EDU/users/aldous/book.html}{http://www.stat.Berkeley.EDU/users/}
  \href{http://www.stat.Berkeley.EDU/users/aldous/book.html}{aldous/book.html}}}.
\bptok{imsref}%
\end{bmisc}
\endbibitem

\bibitem{MR2135733}
\begin{barticle}[mr]
\bauthor{\bsnm{Benjamini},~\bfnm{Itai}\binits{I.}},
  \bauthor{\bsnm{Berger},~\bfnm{Noam}\binits{N.}},
  \bauthor{\bsnm{Hoffman},~\bfnm{Christopher}\binits{C.}} \AND
  \bauthor{\bsnm{Mossel},~\bfnm{Elchanan}\binits{E.}}
(\byear{2005}).
\btitle{Mixing times of the biased card shuffling and the asymmetric exclusion
  process}.
\bjournal{Trans. Amer. Math. Soc.}
\bvolume{357}
\bpages{3013--3029 (electronic)}.
\bid{doi={10.1090/S0002-9947-05-03610-X}, issn={0002-9947}, mr={2135733}}
\bptok{imsref}%
\end{barticle}
\endbibitem

\bibitem{MR2515797}
\begin{barticle}[mr]
\bauthor{\bsnm{Boyd},~\bfnm{Stephen}\binits{S.}},
  \bauthor{\bsnm{Diaconis},~\bfnm{Persi}\binits{P.}},
  \bauthor{\bsnm{Parrilo},~\bfnm{Pablo}\binits{P.}} \AND
  \bauthor{\bsnm{Xiao},~\bfnm{Lin}\binits{L.}}
(\byear{2009}).
\btitle{Fastest mixing {M}arkov chain on graphs with symmetries}.
\bjournal{SIAM J. Optim.}
\bvolume{20}
\bpages{792--819}.
\bid{doi={10.1137/070689413}, issn={1052-6234}, mr={2515797}}
\bptok{imsref}%
\end{barticle}
\endbibitem

\bibitem{MR2202924}
\begin{barticle}[mr]
\bauthor{\bsnm{Boyd},~\bfnm{Stephen}\binits{S.}},
  \bauthor{\bsnm{Diaconis},~\bfnm{Persi}\binits{P.}},
  \bauthor{\bsnm{Sun},~\bfnm{Jun}\binits{J.}} \AND
  \bauthor{\bsnm{Xiao},~\bfnm{Lin}\binits{L.}}
(\byear{2006}).
\btitle{Fastest mixing {M}arkov chain on a path}.
\bjournal{Amer. Math. Monthly}
\bvolume{113}
\bpages{70--74}.
\bid{doi={10.2307/27641840}, issn={0002-9890}, mr={2202924}}
\bptok{imsref}%
\end{barticle}
\endbibitem

\bibitem{MR2124681}
\begin{barticle}[mr]
\bauthor{\bsnm{Boyd},~\bfnm{Stephen}\binits{S.}},
  \bauthor{\bsnm{Diaconis},~\bfnm{Persi}\binits{P.}} \AND
  \bauthor{\bsnm{Xiao},~\bfnm{Lin}\binits{L.}}
(\byear{2004}).
\btitle{Fastest mixing {M}arkov chain on a graph}.
\bjournal{SIAM Rev.}
\bvolume{46}
\bpages{667--689}.
\bid{doi={10.1137/S0036144503423264}, issn={0036-1445}, mr={2124681}}
\bptok{imsref}%
\end{barticle}
\endbibitem

\bibitem{MR2629990}
\begin{barticle}[mr]
\bauthor{\bsnm{Caputo},~\bfnm{Pietro}\binits{P.}},
  \bauthor{\bsnm{Liggett},~\bfnm{Thomas~M.}\binits{T.~M.}} \AND
  \bauthor{\bsnm{Richthammer},~\bfnm{Thomas}\binits{T.}}
(\byear{2010}).
\btitle{Proof of {A}ldous' spectral gap conjecture}.
\bjournal{J. Amer. Math. Soc.}
\bvolume{23}
\bpages{831--851}.
\bid{doi={10.1090/S0894-0347-10-00659-4}, issn={0894-0347}, mr={2629990}}
\bptok{imsref}%
\end{barticle}
\endbibitem

\bibitem{DF}
\begin{barticle}[mr]
\bauthor{\bsnm{Diaconis},~\bfnm{Persi}\binits{P.}} \AND
  \bauthor{\bsnm{Fill},~\bfnm{James~Allen}\binits{J.~A.}}
(\byear{1990}).
\btitle{Strong stationary times via a new form of duality}.
\bjournal{Ann. Probab.}
\bvolume{18}
\bpages{1483--1522}.
\bid{issn={0091-1798}, mr={1071805}}
\bptok{imsref}%
\end{barticle}
\endbibitem

\bibitem{MR94i60074}
\begin{barticle}[mr]
\bauthor{\bsnm{Diaconis},~\bfnm{Persi}\binits{P.}} \AND
  \bauthor{\bsnm{Saloff-Coste},~\bfnm{Laurent}\binits{L.}}
(\byear{1993}).
\btitle{Comparison theorems for reversible {M}arkov chains}.
\bjournal{Ann. Appl. Probab.}
\bvolume{3}
\bpages{696--730}.
\bid{issn={1050-5164}, mr={1233621}}
\bptok{imsref}%
\end{barticle}
\endbibitem

\bibitem{Diekmann97engineeringdiffusive}
\begin{bincollection}[auto:STB|2012/12/05|11:57:16]
\bauthor{\bsnm{Diekmann},~\bfnm{Ralf}\binits{R.}},
  \bauthor{\bsnm{Muthukrishnan},~\bfnm{S.}\binits{S.}} \AND
  \bauthor{\bsnm{Nayakkankuppam},~\bfnm{Madhu~V.}\binits{M.~V.}}
(\byear{1997}).
\btitle{Engineering diffusive load balancing algorithms using experiments}.
In \bbooktitle{Solving Irregularly Structured Problems in Parallel}.
\bseries{Lecture Notes in Computer Science}
\bvolume{1253}
\bpages{111--122}.
\bpublisher{Springer}, \blocation{New York}.
\bptok{imsref}%
\end{bincollection}
\endbibitem

\bibitem{MR373604}
\begin{bbook}[mr]
\bauthor{\bsnm{Feller},~\bfnm{William}\binits{W.}}
(\byear{1968}).
\btitle{An Introduction to Probability Theory and Its Applications. {V}ol.
  {I}},
\bedition{3rd} ed.
\bpublisher{Wiley}, \blocation{New York}.
\bid{mr={0228020}}
\bptok{imsref}%
\end{bbook}
\endbibitem

\bibitem{MR1034418}
\begin{barticle}[mr]
\bauthor{\bsnm{Fiedler},~\bfnm{Miroslav}\binits{M.}}
(\byear{1990}).
\btitle{Absolute algebraic connectivity of trees}.
\bjournal{Linear Multilinear Algebra}
\bvolume{26}
\bpages{85--106}.
\bid{doi={10.1080/03081089008817967}, issn={0308-1087}, mr={1034418}}
\bptok{imsref}%
\end{barticle}
\endbibitem

\bibitem{coarseness}
\begin{barticle}[mr]
\bauthor{\bsnm{Fill},~\bfnm{James~Allen}\binits{J.~A.}}
(\byear{1988}).
\btitle{Bounds on the coarseness of random sums}.
\bjournal{Ann. Probab.}
\bvolume{16}
\bpages{1644--1664}.
\bid{issn={0091-1798}, mr={0958208}}
\bptok{imsref}%
\end{barticle}
\endbibitem

\bibitem{MR1183166}
\begin{barticle}[mr]
\bauthor{\bsnm{Fill},~\bfnm{James~Allen}\binits{J.~A.}}
(\byear{1991}).
\btitle{Time to stationarity for a continuous-time {M}arkov chain}.
\bjournal{Probab. Engrg. Inform. Sci.}
\bvolume{5}
\bpages{61--76}.
\bid{doi={10.1017/S0269964800001893}, issn={0269-9648}, mr={1183166}}
\bptok{imsref}%
\end{barticle}
\endbibitem

\bibitem{MR99g60113}
\begin{barticle}[mr]
\bauthor{\bsnm{Fill},~\bfnm{James~Allen}\binits{J.~A.}}
(\byear{1998}).
\btitle{An interruptible algorithm for perfect sampling via {M}arkov chains}.
\bjournal{Ann. Appl. Probab.}
\bvolume{8}
\bpages{131--162}.
\bid{doi={10.1214/aoap/1027961037}, issn={1050-5164}, mr={1620346}}
\bptok{imsref}%
\end{barticle}
\endbibitem

\bibitem{MR2001m60164}
\begin{barticle}[mr]
\bauthor{\bsnm{Fill},~\bfnm{James~Allen}\binits{J.~A.}},
  \bauthor{\bsnm{Machida},~\bfnm{Motoya}\binits{M.}},
  \bauthor{\bsnm{Murdoch},~\bfnm{Duncan~J.}\binits{D.~J.}} \AND
  \bauthor{\bsnm{Rosenthal},~\bfnm{Jeffrey~S.}\binits{J.~S.}}
(\byear{2000}).
\btitle{Extension of {F}ill's perfect rejection sampling algorithm to general
  chains}.
\bjournal{Random Structures Algorithms}
\bvolume{17}
\bpages{290--316}.
\bptok{imsref}%
\end{barticle}
\endbibitem

\bibitem{holroyddelay}
\begin{barticle}[mr]
\bauthor{\bsnm{Holroyd},~\bfnm{Alexander~E.}\binits{A.~E.}}
(\byear{2011}).
\btitle{Some circumstances where extra updates can delay mixing}.
\bjournal{J. Stat. Phys.}
\bvolume{145}
\bpages{1649--1652}.
\bid{doi={10.1007/s10955-011-0365-x}, issn={0022-4715}, mr={2863724}}
\bptok{imsref}%
\end{barticle}
\endbibitem

\bibitem{MR0356197}
\begin{bbook}[mr]
\bauthor{\bsnm{Karlin},~\bfnm{Samuel}\binits{S.}} \AND
  \bauthor{\bsnm{Taylor},~\bfnm{Howard~M.}\binits{H.~M.}}
(\byear{1975}).
\btitle{A First Course in Stochastic Processes},
\bedition{2nd} ed.
\bpublisher{Academic Press}, \blocation{New York}.
\bid{mr={0356197}}
\bptok{imsref}%
\end{bbook}
\endbibitem

\bibitem{LM}
\begin{barticle}[auto:STB|2012/12/05|11:57:16]
\bauthor{\bsnm{Lange},~\bfnm{Lester~H.}\binits{L.~H.}} \AND
  \bauthor{\bsnm{Miller},~\bfnm{James~W.}\binits{J.~W.}}
(\byear{1992}).
\btitle{A random ladder game: Permutations, eigenvalues, and convergence of
  Markov chains}.
\bjournal{College Math. J.}
\bvolume{23}
\bpages{373--385}.
\bptok{imsref}%
\end{barticle}
\endbibitem

\bibitem{LovaszWinkler2}
\begin{bincollection}[mr]
\bauthor{\bsnm{Lov{\'a}sz},~\bfnm{L{\'a}szl{\'o}}\binits{L.}} \AND
  \bauthor{\bsnm{Winkler},~\bfnm{Peter}\binits{P.}}
(\byear{1995}).
\btitle{Mixing of random walks and other diffusions on a graph}.
In \bbooktitle{Surveys in Combinatorics, 1995 ({S}tirling)}.
\bseries{London Mathematical Society Lecture Note Series}
\bvolume{218}
\bpages{119--154}.
\bpublisher{Cambridge Univ. Press}, \blocation{Cambridge}.
\bid{doi={10.1017/CBO9780511662096.007}, mr={1358634}}
\bptok{imsref}%
\end{bincollection}
\endbibitem

\bibitem{MR81b00002}
\begin{bbook}[mr]
\bauthor{\bsnm{Marshall},~\bfnm{Albert~W.}\binits{A.~W.}} \AND
  \bauthor{\bsnm{Olkin},~\bfnm{Ingram}\binits{I.}}
(\byear{1979}).
\btitle{Inequalities: Theory of Majorization and Its Applications}.
\bseries{Mathematics in Science and Engineering}
\bvolume{143}.
\bpublisher{Academic Press},
  \blocation{New York}.
\bid{mr={0552278}}
\bptok{imsref}%
\end{bbook}
\endbibitem

\bibitem{PeresWinkler}
\begin{bmisc}[auto:STB|2012/12/05|11:57:16]
\bauthor{\bsnm{Peres},~\bfnm{Yuval}\binits{Y.}}
(\byear{2005}).
\bhowpublished{Mixing for Markov chains and spin systems. Lecture notes, Univ.
  British Columbia. Summary available at
  \texttt{\href{http://www.stat.berkeley.edu/\textasciitilde peres/ubc.pdf}{http://www.stat.berkeley.edu/}
  \href{http://www.stat.berkeley.edu/\textasciitilde peres/ubc.pdf}{\textasciitilde peres/ubc.pdf}}}.
\bptok{imsref}%
\end{bmisc}
\endbibitem

\bibitem{PW2011}
\begin{bmisc}[auto:STB|2012/12/05|11:57:16]
\bauthor{\bsnm{Peres},~\bfnm{Yuval}\binits{Y.}} \AND
  \bauthor{\bsnm{Winkler},~\bfnm{Peter}\binits{P.}}
(\byear{2011}).
\bhowpublished{Can extra updates delay mixing? Preprint. Available at
  arXiv:\arxivurl{1112.0603v1} [math.PR]}.
\bptok{imsref}%
\end{bmisc}
\endbibitem

\bibitem{MR99d60081}
\begin{bincollection}[mr]
\bauthor{\bsnm{Propp},~\bfnm{James}\binits{J.}} \AND
  \bauthor{\bsnm{Wilson},~\bfnm{David}\binits{D.}}
(\byear{1998}).
\btitle{Coupling from the past: A user's guide}.
In \bbooktitle{Microsurveys in Discrete Probability ({P}rinceton, {NJ}, 1997)}.
\bseries{DIMACS Series in Discrete Mathematics and Theoretical Computer Science}
\bvolume{41}
\bpages{181--192}.
\bpublisher{Amer. Math. Soc.}, \blocation{Providence, RI}.
\bid{mr={1630414}}
\bptok{imsref}%
\end{bincollection}
\endbibitem

\bibitem{MR99k60176}
\begin{barticle}[mr]
\bauthor{\bsnm{Propp},~\bfnm{James~Gary}\binits{J.~G.}} \AND
  \bauthor{\bsnm{Wilson},~\bfnm{David~Bruce}\binits{D.~B.}}
(\byear{1996}).
\btitle{Exact sampling with coupled {M}arkov chains and applications to
  statistical mechanics}.
\bjournal{Random Structures Algorithms}
\bvolume{9}
\bpages{223--252}.
\bptok{imsref}%
\end{barticle}
\endbibitem

\bibitem{MR99g60116}
\begin{barticle}[mr]
\bauthor{\bsnm{Propp},~\bfnm{James~Gary}\binits{J.~G.}} \AND
  \bauthor{\bsnm{Wilson},~\bfnm{David~Bruce}\binits{D.~B.}}
(\byear{1998}).
\btitle{How to get a perfectly random sample from a generic {M}arkov chain and
  generate a random spanning tree of a directed graph}.
\bjournal{J. Algorithms}
\bvolume{27}
\bpages{170--217}.
\bid{doi={10.1006/jagm.1997.0917}, issn={0196-6774}, mr={1622393}}
\bptok{imsref}%
\end{barticle}
\endbibitem

\bibitem{MR2198603}
\begin{barticle}[mr]
\bauthor{\bsnm{Roch},~\bfnm{S{\'e}bastien}\binits{S.}}
(\byear{2005}).
\btitle{Bounding fastest mixing}.
\bjournal{Electron. Commun. Probab.}
\bvolume{10}
\bpages{282--296 (electronic)}.
\bid{doi={10.1214/ECP.v10-1169}, issn={1083-589X}, mr={2198603}}
\bptok{imsref}%
\end{barticle}
\endbibitem

\bibitem{MR2340874}
\begin{barticle}[mr]
\bauthor{\bsnm{Saloff-Coste},~\bfnm{L.}\binits{L.}} \AND
  \bauthor{\bsnm{Z{\'u}{\~n}iga},~\bfnm{J.}\binits{J.}}
(\byear{2007}).
\btitle{Convergence of some time inhomogeneous {M}arkov chains via spectral
  techniques}.
\bjournal{Stochastic Process. Appl.}
\bvolume{117}
\bpages{961--979}.
\bid{doi={10.1016/j.spa.2006.11.004}, issn={0304-4149}, mr={2340874}}
\bptok{imsref}%
\end{barticle}
\endbibitem

\bibitem{MR2519527}
\begin{barticle}[mr]
\bauthor{\bsnm{Saloff-Coste},~\bfnm{L.}\binits{L.}} \AND
  \bauthor{\bsnm{Z{\'u}{\~n}iga},~\bfnm{J.}\binits{J.}}
(\byear{2009}).
\btitle{Merging for time inhomogeneous finite {M}arkov chains. {I}. {S}ingular
  values and stability}.
\bjournal{Electron. J. Probab.}
\bvolume{14}
\bpages{1456--1494}.
\bid{doi={10.1214/EJP.v14-656}, issn={1083-6489}, mr={2519527}}
\bptok{imsref}%
\end{barticle}
\endbibitem

\bibitem{SCZwave}
\begin{barticle}[mr]
\bauthor{\bsnm{Saloff-Coste},~\bfnm{L.}\binits{L.}} \AND
  \bauthor{\bsnm{Z{\'u}{\~n}iga},~\bfnm{J.}\binits{J.}}
(\byear{2010}).
\btitle{Time inhomogeneous {M}arkov chains with wave-like behavior}.
\bjournal{Ann. Appl. Probab.}
\bvolume{20}
\bpages{1831--1853}.
\bid{doi={10.1214/09-AAP661}, issn={1050-5164}, mr={2724422}}
\bptok{imsref}%
\end{barticle}
\endbibitem

\bibitem{SCZmergingII}
\begin{barticle}[mr]
\bauthor{\bsnm{Saloff-Coste},~\bfnm{L.}\binits{L.}} \AND
  \bauthor{\bsnm{Z{\'u}{\~n}iga},~\bfnm{J.}\binits{J.}}
(\byear{2011}).
\btitle{Merging for inhomogeneous finite {M}arkov chains, {P}art {II}: {N}ash
  and log-{S}obolev inequalities}.
\bjournal{Ann. Probab.}
\bvolume{39}
\bpages{1161--1203}.
\bid{doi={10.1214/10-AOP572}, issn={0091-1798}, mr={2789587}}
\bptnote{check year}%
\bptok{imsref}%
\end{barticle}
\endbibitem

\bibitem{MR2278445}
\begin{barticle}[mr]
\bauthor{\bsnm{Sun},~\bfnm{Jun}\binits{J.}},
  \bauthor{\bsnm{Boyd},~\bfnm{Stephen}\binits{S.}},
  \bauthor{\bsnm{Xiao},~\bfnm{Lin}\binits{L.}} \AND
  \bauthor{\bsnm{Diaconis},~\bfnm{Persi}\binits{P.}}
(\byear{2006}).
\btitle{The fastest mixing {M}arkov process on a graph and a connection to a
  maximum variance unfolding problem}.
\bjournal{SIAM Rev.}
\bvolume{48}
\bpages{681--699}.
\bid{doi={10.1137/S0036144504443821}, issn={0036-1445}, mr={2278445}}
\bptok{imsref}%
\end{barticle}
\endbibitem

\bibitem{MR1630416}
\begin{bincollection}[mr]
\bauthor{\bsnm{Wilson},~\bfnm{David~B.}\binits{D.~B.}}
(\byear{1998}).
\btitle{Annotated bibliography of perfectly random sampling with {M}arkov
  chains}.
In \bbooktitle{Microsurveys in Discrete Probability ({P}rinceton, {NJ}, 1997)}.
\bseries{DIMACS Series in Discrete Mathematics and Theoretical Computer Science}
\bvolume{41}
\bpages{209--220}.
\bpublisher{Amer. Math. Soc.}, \blocation{Providence, RI}.
\bnote{Latest updated version is posted at \url{http://dbwilson.com/exact/}.}
\bid{mr={1630416}}
\bptok{imsref}%
\end{bincollection}
\endbibitem

\bibitem{MR2001h65014}
\begin{bincollection}[mr]
\bauthor{\bsnm{Wilson},~\bfnm{David~Bruce}\binits{D.~B.}}
(\byear{2000}).
\btitle{Layered multishift coupling for use in perfect sampling algorithms
  (with a primer on {CFTP})}.
In \bbooktitle{Monte {C}arlo Methods ({T}oronto, {ON}, 1998)}.
\bseries{Fields Institute Communications}
\bvolume{26}
\bpages{143--179}.
\bpublisher{Amer. Math. Soc.}, \blocation{Providence, RI}.
\bid{mr={1772313}}
\bptok{imsref}%
\end{bincollection}
\endbibitem

\end{thebibliography}
\end{document}